\documentclass{article}
\usepackage{amsmath}
\usepackage{amssymb}
\usepackage{color}

\usepackage{latexsym}
\usepackage[english]{babel}
\usepackage[dvips]{graphicx}
\usepackage{verbatim}
\newcommand{\hhh}{}

\newcommand{\bbibitem}{\bibitem}

\newcommand{\llabel}[1]{{\label{#1}}}

\newcommand{\ffoot}[1]{}

\renewcommand{\r}[1]{(\ref{#1})}

\newcommand{\ex}[1]{}

\font\tenmsb=msbm10
\font\sevenmsb=msbm7
\font\fivemsb=msbm5
\newfam\msbfam
\textfont\msbfam=\tenmsb
\scriptfont\msbfam=\sevenmsb
\scriptscriptfont\msbfam=\fivemsb
\def\Bbb#1{{\fam\msbfam\relax#1}}

\newcommand{\bi}{\begin{itemize}}
\newcommand{\ei}{\end{itemize}}
\newcommand{\bd}{\begin{description}}
\newcommand{\ed}{\end{description}}
\renewcommand{\i}{\item}

\newcommand{\bqn}{\begin{eqnarray}}
\newcommand{\eqn}{\end{eqnarray}}
\newcommand{\eqnn}{\nonumber\end{eqnarray}}
\newcommand{\eqnl}[1]{\llabel{#1}\end{eqnarray}}

\newcommand{\nn}{\nonumber}
\newcommand{\ba}[1]{\begin{array}{#1}}
\newcommand{\ea}{\end{array}}

\newcommand{\R}{\Bbb{R}}
\newcommand{\mc}{\Bbb{C}}

\newcommand{\fine}{\end{document}}

\def \trait (#1) (#2) (#3){\vrule width #1pt height #2pt depth #3pt}
\def \qed{\hfill
        \trait (0.1) (6) (0)
        \trait (6) (0.1) (0)
        \kern-6pt
        \trait (6) (6) (-5.9)
        \trait (0.1) (6) (0)
\medskip}
\def \qedmio{\hfill
             \trait (8) (8) (-0.1)
             \medskip}
\def \quadp{{\Huge $\qedmio$}}

\newtheorem{ml}{\bf Lemma}
\newtheorem{Theorem}{\bf Theorem}
\newtheorem{mo}{\bf \underline{{\sl Observation}}}
\newtheorem{mcc}{\bf Corollary}
\newtheorem{Definition}{\bf Definition}
\newtheorem{mpr}{\bf Proposition}
\newtheorem{mproperty}{\bf Property}

\newcommand{\bt}{\begin{Theorem}}
\newcommand{\et}{\end{Theorem}}
\newcommand{\bl}{\begin{ml}}
\newcommand{\el}{\end{ml}}
\newcommand{\bo}{\noindent\begin{mo}\rm}
\newcommand{\eo}{\end{mo}}
\newcommand{\bp}{\begin{mpr}}
\newcommand{\ep}{\end{mpr}}
\newcommand{\bc}{\begin{mcc}}
\newcommand{\ec}{\end{mcc}}
\newcommand{\bdeff}{\begin{Definition}}
\newcommand{\edeff}{\end{Definition}}
\newcommand{\bproperty}{\begin{mproperty}}
\newcommand{\eproperty}{\end{mproperty}}

\newtheorem{mrem}{\bf \underline{{\sl Remark}}}
\newcommand{\brem}{\begin{mrem}\rm}
\newcommand{\erem}{\end{mrem}}

\newcommand{\ppotR}[3]
{

\begin{figure}\begin{center}
~\includegraphics[width=#3truecm]{./#1.eps}\\
\caption{#2}
\llabel{#1}
\end{center}
\end{figure}
\noindent$\!\!$}

\newcommand{\lam}{\lambda}
\newcommand{\g}{\gamma}
\newcommand{\al}{\alpha}
\newcommand{\eps}{\varepsilon}
\newcommand{\de}{\delta}

\newcommand{\A}{{\cal A}}

\newcommand{\U}{{\cal U}}
\newcommand{\Supp}{\mbox{{\rm Supp}}}
\newcommand{\ga}{\gamma}

\newcommand{\qq}{{\cal Z}}

\newcommand{\ph}{\varphi}
\renewcommand{\c}{{\cos(\ph)}}

\newcommand{\proof}{{\bf Proof. }}

\newcommand{\para}{\frac{2}{\sqrt{|\delta|}   }}
\newcommand{\paraa}{\frac{2}{\sqrt{\delta}   }}

\newcommand{\rp}{\R P^1}

\begin{document}
\begin{center} \noindent
{\LARGE{\sl{\bf Stability of Planar Switched Systems: the
Nondiagonalizable Case}}}
\end{center}

\vskip 1cm
\begin{center}
Moussa BALDE,

{\footnotesize Dept. de Mathé\'ematiques et Informatique,
Universit\'e Cheikh Anta Diop, Dakar-Fann, Senegal}

\vskip .2cm
Ugo BOSCAIN,

{\footnotesize SISSA-ISAS, Via Beirut 2-4, 34014 Trieste, Italy}

\end{center}

\vspace{2cm} \noindent \rm
\begin{quotation}
\noindent  {\bf Abstract}
%%%%%%%%%%%%%%%%%%%%%%%%%%%%%%%%%%%%%%%%%%%%%%%%%%%%%%%%%%%%%%%%%%%%%%%%%%%%
Consider the planar linear switched system
$\dot x(t)=u(t)Ax(t)+(1-u(t))Bx(t),$ where $A$ and $B$ are two $2\times2$
real matrices, $x\in~\R^2$, and $u(.):[0,\infty[\to\{0,1\}$ is a
measurable function. In this paper
we
consider the problem of finding a (coordinate-invariant) necessary and
sufficient condition on $A$ and $B$
under which the system is asymptotically  stable for arbitrary
switching functions $u(.)$.

This problem was solved in previous works under the assumption
that both $A$ and $B$ are diagonalizable.  In this paper we conclude this
study, by providing a necessary and sufficient condition for asymptotic
stability in the case in which $A$ and/or $B$ are not diagonalizable.

To this purpose we build suitable normal forms for $A$ and $B$ containing
coordinate invariant parameters. A necessary and sufficient condition is
then found without looking for a common  Lyapunov function but using
``worst-trajectory'' type arguments.
%%%%%%%%%%%%%%%%%%%%%%%%%%%%%%%%%%%%%%%%%%%%%%%%%%%%%%%%%%%%%%%%%%%%%%%%%%%%
\end{quotation}

{\bf AMS: }
93D20, 37N35

\begin{center}
PREPRINT SISSA 44/2006/M
\end{center}

\newpage
\section{Introduction}
%In recent years, the problem of stability and stabilizability of switched
%systems has attracted increasing attentions
%(see for instance
%\cite{agr,sw-1,DM,liberzon-book,survey}).

By a  switched system, we mean a family of continuous--time dynamical
systems and a rule that determines at each time which dynamical system
 is
responsible of the time evolution. More precisely,
let $\{f_u:~u\in U\}$ (where $U$ is a subset of
$\R^m$) be a finite or infinite
set of sufficiently regular vector fields  on
a manifold $M$, and consider the family of dynamical systems:
\bqn
\dot x=f_u(x),~~x\in M.
\eqnl{sw-nl}
The rule  is given by assigning the so-called switching
function, i.e., a function $u(.):[0,\infty[\to U\subset
\R^m$. Here, we consider the situation in which
the switching function is not known a priori and represents some
phenomenon (e.g., a disturbance) that is not possible to control.
Therefore, the dynamics defined in \r{sw-nl} also fits into the framework
of uncertain systems (cf. for instance \cite{bub}).
In the sequel, we use the notations $u\in U$ to label a fixed
individual system and $u(.)$ to indicate the switching function.
These kind of systems are sometimes called ``n-modal systems'',
``dynamical poly\-systems'',  ``input systems''.
The term ``switched system'' is often reserved to situations in which
{\hhh the switching function $u(.)$ is piecewise continuous or} the set
$U$
is finite. For the purpose of this paper, we only
require $u(.)$ to be a measurable function.
For a discussion of various issues related
to switched systems, we refer the reader to
\cite{blanchini,DM,liberzon-book,survey}.

A typical problem for switched systems goes as follows.
Assume that, for every fixed $u\in U$, the dynamical system $\dot
x=f_u(x)$
satisfies
a given property (P). Then one can investigate conditions under which
property (P) still holds for $\dot x=f_{u(t)}(x)$, where $u(.)$ is an
arbitrary
switching function.

\bigskip
\noindent
In \cite{agr,sw-1,DM,hesp}, the case of linear switched systems was
considered:
\bqn
\label{sw-lup}\dot x(t)=A_{u(t)}
x(t),~~x\in\R^n,~~\{A_u\}_{u\in U}\subset\R^{n\times n},
\eqnl{-2}
where $U\subset\R^m$ is a compact set,
$u(.):[0,\infty[\to U$ is a (measurable) switching function, and the map
$u\mapsto A_u$ is continuous (so that
$\{A_u\}_{u\in U}$ is a compact set of matrices).
For these systems, the problem of asymptotic stability of the origin,
uniformly with respect to switching functions  was
investigated.

Let us recall the notions of stability that are used in the following.
%%%%%%%%%%%%%%%%%%%%%%%%%%%%%%%%%%%%%%%%%%%%%%%%%%%%%%%%%%%%%%%%%%%%%%%%%%%
\bdeff
\label{d-stability}
For $\de>0$ let  $B_\de$ be the unit ball of
radius $\de$, centered in the origin. Denote by
$\U$ the set of measurable
functions defined on $[0,\infty[$  and taking values on $U$. Given
$x_0\in \R^n$, we denote by $\g_{x_0,u(.)}(.)$ the trajectory of \r{-2}
based in
$x_0$ and corresponding to the control $u(.)$.
The {\it accessible set from $x_0$},
denoted by ${\A(x_0)}$, is
$$\A(x_0)=\cup_{u(.)\in\U}\Supp(\ga_{x_0,u(.)}(.))\,.$$
%%%%%%%%%%%%%%%%%%%%%%%%%%%%%%%%%%%%%%%%%%%%%%%%%%%%%%%%%%%
We say that the system \r{-2} is
\bi
\i {\bf unbounded}
   if there exist $x_0\in \R^n$
and $u(.)\in\U$ such that
 $\ga_{x_0,u(.)}(t)$ goes to infinity as $t\to\infty$;
\i {\bf uniformly stable}
if, for every $\eps>0$, there exists $\de>0$ such that
$\A(x_0)\subset B_\eps$ for every
$x_0\in B_\de$;
\i {\bf globally uniformly asymptotically stable} ({\bf GUAS}, for short)
if it is
uniformly
stable and globally uniformly attractive,  i.e.,
for every $\de_1, \de_2>0$, there exists $T>0$ such that
$\ga_{x_0,u(.)}(T)\in
B_{\de_1}$ for every
$u(.)\in\U$ and every $x_0\in B_{\de_2}$;
\ei
\edeff
%%%%%%%%%%%%%%%%%%%%%%%%%%%%%%%%%%%%%%%%%%%%%%%%%%%%%%%%%%%%%%%%%%%
\brem
Under our hypotheses (linearity and compactness) there are many notions of
stability equivalent to the ones of Definition \ref{d-stability}.
More precisely since the system is linear, local and global notions of
stability are equivalent. Moreover, since $\{A_u \}_{u\in U}$ is compact,
all notions of stability are
automatically uniform with respect to switching functions (see for
instance \cite{AISW}).
Finally, thanks to the linearity, the GUAS property is equivalent to
the more often quoted property of GUES (global exponential stability,
uniform with respect to switching), see for example \cite{angeli} and
references therein.
\erem
%%%%%%%%%%%%%%%%%%%%%%%%%%%%%%%%%%%%%%%%%%%%%%%%%%%%%%%%%%%%%%%%%%%%%
Let us recall some results about stability of systems of type \r{-2}.

In \cite{agr,hesp}, it is shown that
the structure of the Lie algebra generated by the
matrices $A_u$:
\bqn
\mbox{{\bf g}}=\{A_u:~u\in U\}_{L.A.},
\eqnn
is crucial for the stability of the system \r{-2}.
The main result of
\cite{hesp} is the following:
%%%%%%%%%%%%%%%%%%%%%%%%%%%%%%%%%%%%%%%%%%%%%%%%%%%%%%
\bt
{\bf (Hespanha, Morse, Liberzon)} If $\mbox{{\bf g}}$ is a solvable Lie
algebra, then
the
switched system \r{-2} is GUAS.
\label{t-hes}
\et

In \cite{agr} a generalization was given. Let $\mbox{{\bf
g}}=\mbox{{\bf r}}~\ltimes~
\mbox{{\bf s}}$ be the Levi decomposition of $\mbox{{\bf g}}$ in its
radical
(i.e.,
the maximal solvable ideal of $\mbox{{\bf g}}$) and a semi--simple
sub--algebra,
where  the symbol $\ltimes$ indicates the semidirect sum.
\bt
{\bf (Agrachev, Liberzon)} If $\mbox{{\bf s}}$ is a compact Lie algebra
then
the
switched system \r{-2} is GUAS.
\label{t-agr}
\et
%%%%%%%%%%%%%%%%%%%%%%%%%%%%%%%%%%%%%%%%%%%%%%%%%%%%%%%%%%%%%
Theorem \ref{t-agr} contains Theorem \ref{t-hes} as a special case. Anyway
the converse of Theorem \ref{t-agr} is not true in general: if
$\mbox{{\bf s}}$
is non compact, the system can be stable or unstable. This case was
also investigated.
In particular, if $\mbox{{\bf g}}$ has dimension at most $4$ as Lie
algebra,
the authors were able to
reduce
the problem of the asymptotic stability of the system \r{-2}
to the problem of
the
asymptotic stability of an auxiliary bidimensional system. We refer the
reader to \cite{agr} for details. For this reason the bidimensional
problem assumes particularly interest. In \cite{sw-1} (see also
\cite{sw-lyapunov}) the single input
case was investigated,
\bqn
\dot x(t)=u(t)Ax(t)+(1-u(t))Bx(t),
\eqnl{1}
where $A$ and $B$ are two $2\times2$ real matrices with eigenvalues having
strictly negative real part (Hurwitz in the following), $x\in\R^2$ and
$u(.):[0,\infty[\to\{0,1\}$
is an arbitrary  measurable switching function.

Under the assumption that $A$ and $B$ are both diagonalizable (in real or
complex sense)
a complete solution was found. In the
following we refer to this case as to the {\it diagonalizable case}.
More precisely a necessary and sufficient condition for GUAS was given
in terms of three coordinate-invariant parameters: one
 depends on the eigenvalues of $A$,
 one
on the eigenvalues of $B$
and the last
contains the
interrelation
among the two systems and it is in 1--1 correspondence with the cross
ratio
of the four
eigenvectors
of $A$ and $B$
in the projective line
$\mc P^1$.
A remarkable fact is that if the system $\dot x=u(t)Ax+(1-u(t))Bx$ has a
given stability property, then  the system $\dot
x=u(t)(\tau_1A)x+(1-u(t))(\tau_2B)x$ has the same stability
property, for every $\tau_1,\tau_2>0$. This is a consequence of the fact
that the stability
properties of the system \r{1}, depend only on the shape of the
integral curves of $Ax$ and $Bx$ and not on the way
in which they are parameterized.

The stability conditions for \r{1} were obtained with a direct method
without looking for a
common Lyapunov function, but analyzing the locus in which the two
vector fields are collinear, to build the ``worst-trajectory''.
This method was successful also to study a nonlinear
generalization of this problem (see \cite{sw-nonlinear}).
One of
the most important step to obtain these stability conditions was
to find good normal forms for the matrices, containing explicitly the
coordinate invariant parameters.

%%%%%%%%%%%%%%%%%%%%%%%%%%%%%%%%%%%%%%%%%%%%%%%%%%%%%%%%%%%%%%%%%%%%%%%%%%%
\brem
It is interesting to notice that common Lyapunov functions do not seem
to be
the most efficient tool to study the stability of switching systems for
arbitrary switchings. In fact, beside providing sufficient conditions for
GUAS, like those of Theorems \ref{t-hes} and \ref{t-agr}, the concept of
common
Lyapunov function is useful when one
can prove that, if a Lyapunov function
exists, then it is possible
to find it in a class of functions parameterized by a finite number
of parameters.
Indeed, once such a class of functions is identified, then in order to
verify
GUAS, one could use numerical algorithms to check (by varying the
parameters) whether a Lyapunov function exists (in which case the
system is GUAS) or not (meaning that the system is not GUAS).
In the case of stability of switching systems under
arbitrary switchings, this is in general a very hard task. For instance
for systems of type \r{1} one can prove that
the GUAS property is equivalent to the existence of a common {\it
polynomial} Lyapunov function,  but the degree of such common
polynomial Lyapunov
function is not uniformly bounded over all
the GUAS systems \cite{sw-lyapunov}.
\erem
%%%%%%%%%%%%%%%%%%%%%%%%%%%%%%%%%%%%%%%%%%%%%%%%%%%%%%%

In this paper we provide a necessary and sufficient condition for
GUAS  in the nongeneric cases omitted in \cite{sw-1}. In
particular we study the stability of the system \r{1} assuming that at least one of
the matrices (say $A$) is not diagonalizable. In the following we refer to
this case as to the {\it nondiagonalizable case}.
We also assume that $A$ and
$B$ are both Hurwitz and that  $[A,B]\neq0$ otherwise the problem is
trivial. These conditions are gathered in the following assumption {\bf
(H0)}, which is often
recalled in the following:
\bd
\i[(H0)] $A$ and $B$ are two $2\times2$ real
Hurwitz matrices.
We assume that $A$ is nondiagonalizable and $[A,B]\neq0$.
\ed
We also study the cases in which the system is just
uniformly stable.\\\\
%%%%%%%%%%%%%%%%%%%%%%%%%%%%%%%%%%
A very useful fact is that the stability properties of systems of kind
\r{-2} depend only on the convex hull of the set $\{A_u\}_{u\in U}$ (see
for instance \cite{sw-lyapunov}). As a
consequence we have
%%%%%%%%%%%%%%%%%%%%%%%%%%%%
\bl The system \r{1} with $u(.):[0,\infty[\to \{0,1\}$ is GUAS (resp.
uniformly stable, resp. unbounded) if and only the system \r{1} with
$u(.):[0,\infty[\to [0,1]$ is.
\el
%%%%%%%%%%%%%%%%%%%%%%%%%%%%
In the following we refer to the switched system with $u(.)$ taking
values in $[0,1]$ as to the {\it convexified system}.
Sometimes we will take advantage of studying the convexified system.\\
%Thanks to this lemma, we can study the stability properties of the
%convexified system when it is useful.\\

The techniques that we use to get the stability conditions for $\r{1}$
under
{\bf (H0)} are similar to those of  \cite{sw-1}. However new difficulties
arise. The
first
is due to the fact that since $A$ is not diagonalizable then the
eigenvectors of $A$ and $B$ are at most 3 noncoinciding points on $\mc
P^1$. As a
consequence the cross ratio is not anymore the right
parameter describing the interrelation among the systems. It is
either not defined or completely fixed. For this reason  new
coordinate-invariant parameters should be identified and new normal
forms for  $A$ and $B$ should be constructed. These coordinate
invariant parameters are the three real parameters defined in Definition
\ref{d-invariants} below.
One ($\eta$) is the (only) eigenvalue of $A$, the second ($\rho$)
depends on the eigenvalues of $B$ and the third ($k$) plays the role of the
cross ratio of the
diagonalizable case. Section \ref{s-normal} is devoted to the
computation of the normal forms.\\

Once suitable normal forms are obtained, we look for stability conditions, studying the set
$\qq$ where the two vector fields are
linearly dependent, using a technique coming from optimal control (see \cite{libro}). This set is the
set of zeros of the function
$Q(x):=\det(Ax,Bx)$. Since $Q$ is a quadratic form, we have the following
cases (depicted in Figure \ref{f-stabbb}):
\bd
\i[A.] $\qq=\{0\}$
(i.e., $Q$ is positive or negative definite).
In this case one
vector field points always on the
same
side of the other and the system is GUAS. This fact can be proved in
several way (for instance building a common quadratic Lyapunov
function)  and it is true in much more generality (even for nonlinear
systems, see \cite{sw-nonlinear}).
\i[B.] $\qq$ is the union of two noncoinciding straight lines passing
through
the origin
(i.e., $Q$ is sign indefinite).
Take a point $x\in\qq\setminus\{0\}$.
We say that $\qq$ is {\it direct} (respectively, {\it inverse}) if $Ax$
and $Bx$
have the same (respectively, opposite) versus.
One can prove that this definition is independent of the choice of $x$ on
$\qq$.
See Proposition \ref{p-or} below. Then we have the two subcases:
\bd
\i[B1.] $\qq$ is inverse.
In this
case one can prove that there exists $u_0\in]0,1[$ such that the matrix
$u_0Ax+(1-u_0)Bx$ has an eigenvalue with positive real part. In this
case the system is unbounded since it is possible to build a
trajectory of the convexified system going to infinity with constant
control. (This
type of instability is called {\it static instability}.)

\i[B2.] $\qq$ is direct. In this
case one can reduce the problem of the stability of \r{1} to the problem
of
the stability of a single trajectory called {\it worst-trajectory}.
Fixed $x_0\in\R^2\setminus\{0\}$, the worst-trajectory
$\g_{x_0}$ is the trajectory of \r{1}, based at
$x_0$, and  having the following
property. At each time $t$, $\dot \g_{x_0}(t)$ forms the
smallest angle (in absolute value) with the
(exiting) radial direction (see Figure \ref{f-worst}). Clearly the
worst-trajectory
switches among the two vector fields on the set $\qq$.
If it does not rotate around the origin
(i.e., if it crosses the set $\qq$ a finite number of times) then the
system is GUAS. This case is better described projecting the system on $\rp$ (see Lemma
\ref{l-projective} below).
On the other
side,  if it
rotates around the origin, the system is GUAS if and only if after
one turn the distance from the origin is decreased.
(see Figure \ref{f-stabbb}, Case B2). If after one turn the
distance from the origin is increased then the system is unbounded (in
this case, since there are no trajectories of the convexified system
going to
infinity with constant control, we call this instability  {\it dynamic instability}).
If  $\g_{x_0}$ is periodic then the system is uniformly stable, but not
GUAS.
\ed
\i[C.] In the degenerate case in which the two straight lines of $\qq$
coincide (i.e., when $Q$ is sign semi-definite),
one see that the system is GUAS (resp. uniformly stable, but
not GUAS) if and
only if $\qq$ is direct (resp. inverse). We call these cases respectively
{\bf C2} and {\bf C1}.
\ed
%%%%%%%%%%%%%%%%%%%%%%%%%%%%
\ppotR{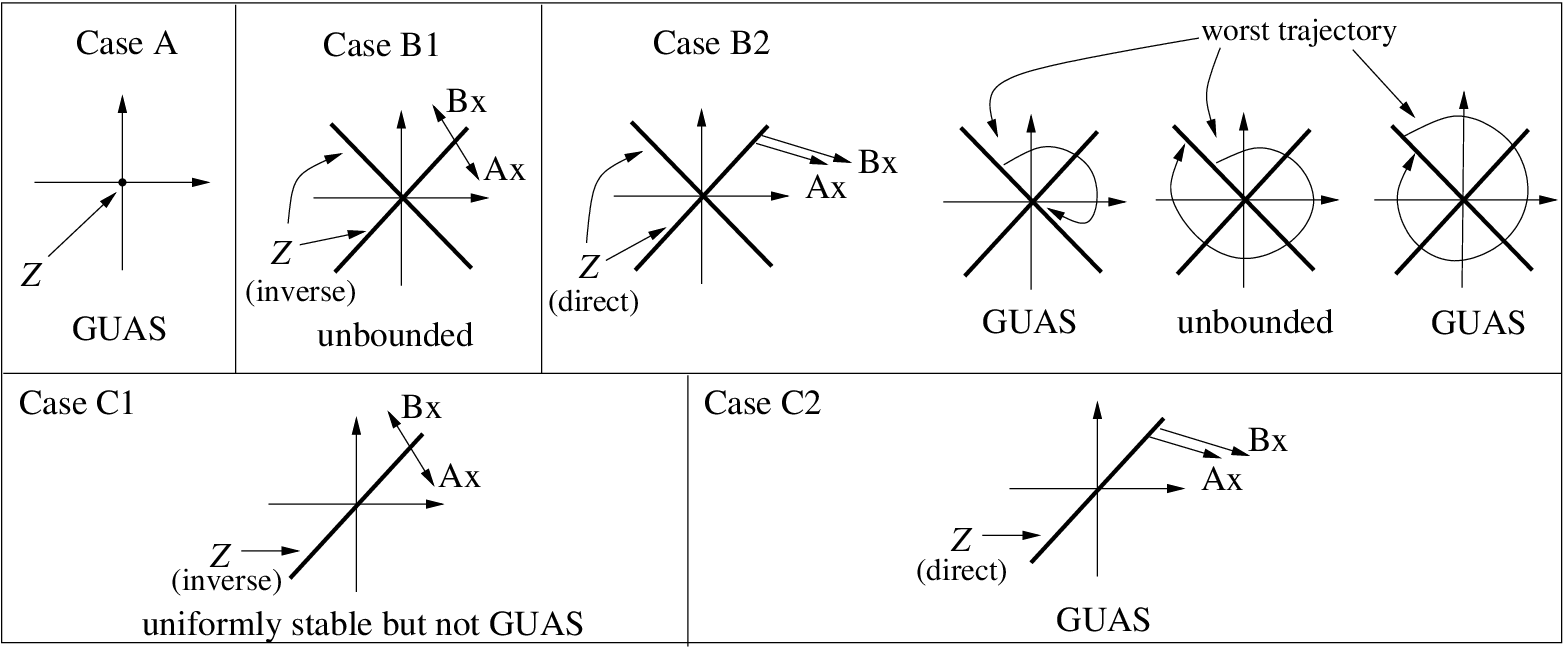}{}{16}
%%%%%%%%%%%%%%%%%%%%%%%%%%%%
The main point of the paper is to translate conditions
{\bf A.},
{\bf B.} and
{\bf C.} in terms of the coordinate invariant parameters appearing in the normal forms. Several
cases (parametrized by the parameter $k$) should be studied separately:
\bi
\i $k=0$, called singular case ({\bf S}-case in the following),
\i $k\neq0$ and  $B$ has non real eigenvalues, called ``regular~$-1$'' case ({\bf R}$_{-1}$ case),
\i $k\neq0$ and  $B$ has real noncoinciding eigenvalues
called ``regular~$+1$'' case ({\bf R}$_{1}$ case),
\i $k\neq0$ and  $B$ is not diagonalizable, called ``regular~$0$'' case ({\bf R}$_{0}$ case).
\ei
The case in which $B$ has two coinciding
eigenvalues and it is diagonalizable is not considered since in
this case $B$ is proportional to the identity and, therefore, {\bf (H0)}
is
not satisfied.\\

The structure of the paper is the following.
In Section \ref{s-normal} we compute the normal forms for $A$ and $B$ in
which the coordinate-invariant parameters appear explicitly.
In Section \ref{s-main} we state the stability conditions, that are proved
in the next sections.
%%%%%%%%%%%%%%%%%%%%%%%
\brem
As in the diagonalizable case (since the way in which the integral curves
of $Ax$ and $Bx$ are parametrized is not important), if the system $\dot
x=u(t)Ax+(1-u(t))Bx$
has a given stability property, then  the system $\dot
x=u(t)(\tau_1A)x+(1-u(t))(\tau_2B)x$ has the same stability
property, for every $\tau_1,\tau_2>0$.
\erem
%%%%%%%%%%%%%%%%%%%%%%%%%%
In section
\ref{s-parallel} we start by
studying the set  $\qq$ where the two vector
fields are linearly dependent.
In Section \ref{s-general} we state and prove some
general stability conditions (in particular cases {\bf A}, {\bf B1}, {\bf C1}, {\bf C2} of Figure
\ref{f-stabbb}).
In Section \ref{s-worst} we build the worst-trajectory (i.e., we study
case
{\bf B2}).
\ppotR{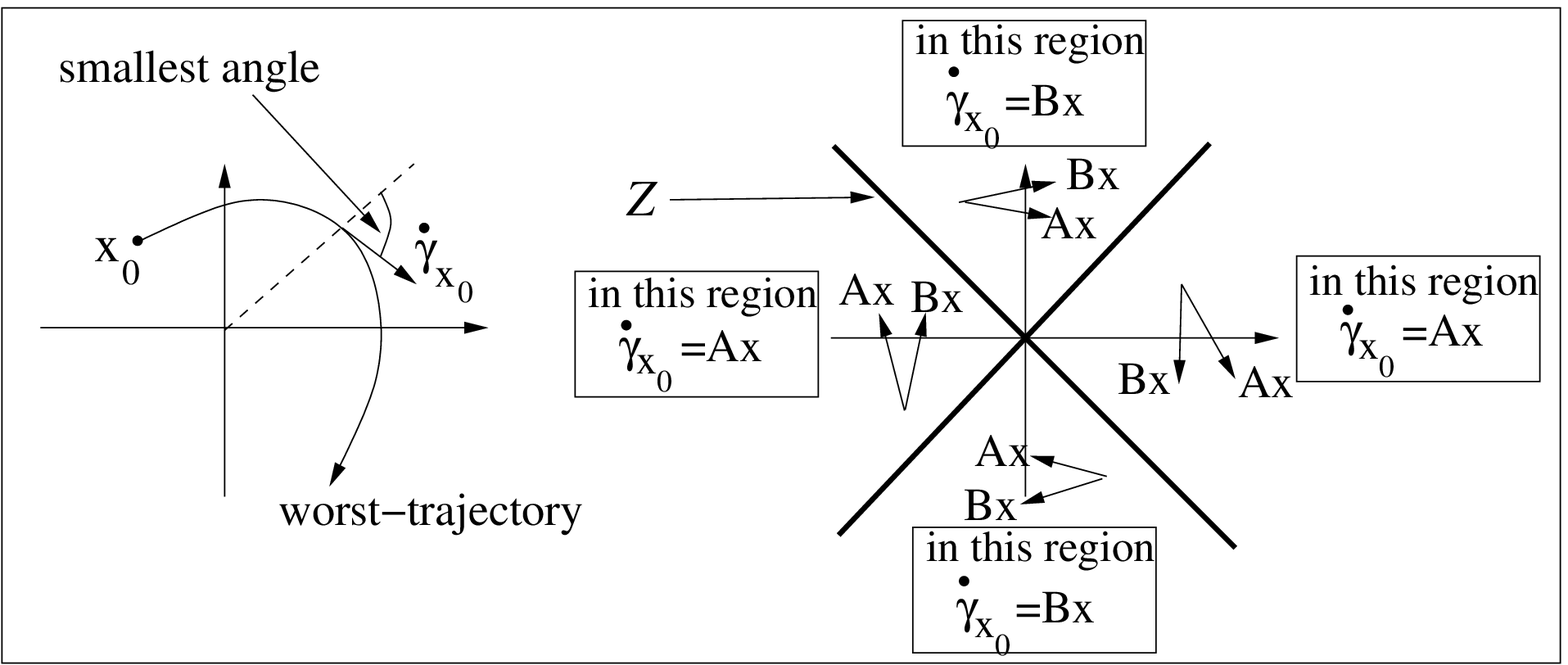}{}{12}

%%%%%%%%%%%%%%%%%%%%%%%%%%%%%%%%%%%%%%%%%%%%
\section{Basic Definitions and Normal Forms}
\llabel{s-normal}
For $x\in\R$ define
$$sign(x)=\left\{\ba{l}+1\mbox{ if $x>0$  }\\
0\mbox{ if $x=0$ }\\
-1\mbox{ if $x<0$.  }
\ea\right.
$$
%%%%%%%%%%%%%%%%%%%%%%%%%%%%%%%%
\bdeff
\llabel{d-invariants}
Assume {\bf (H0)} and let $\delta$ be the discriminant
of the equation $det(B-\lam $Id$)=0$.
Define the following invariant
parameters:
\bqn
&&
\eta=\left\{\ba{l}\displaystyle\frac{Tr(A)}{\sqrt{|\delta|}}\mbox{ if
}\delta\neq0\\
\displaystyle\frac{Tr(A)}{2}\mbox{ if }\delta=0,\ea\right.\\
&&\rho=\left\{\ba{l}\displaystyle\frac{Tr(B)}{\sqrt{|\delta|}}\mbox{ if
}\delta\neq0\label{e-rho}\\
\displaystyle\frac{Tr(B)}{2}\mbox{ if }\delta=0,\ea\right.\\
&&k=\left\{\ba{l}\displaystyle
\frac{4}{|\delta|}\left(Tr(AB)-\frac12Tr(A)
Tr(B)\right)\mbox{ if
}\delta\neq0\\
\displaystyle Tr(AB)-\frac12Tr(A)
Tr(B)\mbox{ if
}\delta=0.\     \ea\right.\label{e-k}\eqn
\edeff
%%%%%%%%%%%%%%%%%%%%%%%%%%%%%%%%%%%%%%%%%%%%%%%%%%%%%%%%%%%%
\brem Notice that $\delta=(\lam_1-\lam_2)^2\in\R$, where $\lam_1$
and $\lam_2$ are the eigenvalues of $B$. Notice moreover that $B$
has non real eigenvalues if and only if $\delta<0$. Finally
observe that $\eta,\rho<0$ and $k\in\R$. \erem
%%%%%%%%%%%%%%%%%%%%%%%%%%%%%%%%%%%%%%%%%%%%%%%%%%%%%%%%%%%%
\bdeff
In the following, under the assumption {\bf (H0)}, we call  regular case
({\bf R}-case for short), the case in which $k\neq0$ and singular case
({\bf S}-for short), the case in which $k=0$.
\edeff

%%%%%%%%%%%%%%%%%%%%%%%%%%%%%%%%%%%%%%%%%%%%%%%%%%%%%%%
\bl {\bf (R-case)}
\label{l-nf-kneq0}
Assume {\bf (H0)} and $k\neq0$.
Then it is always possible to find a linear
change of
coordinates and a constant $\tau>0$ such that $A/\tau$ and $B/\tau$
(that we still call $A$ and $B$) have the following form:
\bqn
&&A=\left(\ba{cc}\eta&1\\0&\eta\ea\right),\label{e-nA}\\
&&B=\left(\ba{cc}\rho&sign(\delta)/k\\k&\rho\ea\right)\label{e-nB}.
\eqn
Moreover in this case $[A,B]\neq0$ is automatically satisfied.
\el
%%%%%%%%%%%%%%%%%%%%%%%%%%%%%%%%%%%%%%%%%%%%%%%%%%%%%%%%
\brem
In the following, for the regular case,
we call
{\bf R}$_{1}$-case,
{\bf R}$_{-1}$-case,
{\bf R}$_{0}$-case the cases corresponding respectively to
$sign(\delta)=1$,
$sign(\delta)=-1$,
$sign(\delta)=0$. See Lemma \ref{l-autovalori} below for the discussion of
the
eigenvalues and eigenvectors of $B$ in these three cases.
\erem
%%%%%%%%%%%%%%%%%%%%%%%%%%%%%%%%%%%%%%%%%%%%%%%%%%%%%%%%%%%%%%%%%%%
\bl {\bf (S-case)}
\label{l-nf-k=0}
Assume {\bf (H0)} and $k=0$. Then $\delta>0$ and
it is always possible to find a linear
change of
coordinates and a constant $\tau>0$ such that $A/\tau$ and $B/\tau$
(that we still call $A$ and $B$) have the following form,
\bqn
A=\left(\ba{cc}\eta&1\\0&\eta\ea\right),~~~
B=\left(\ba{cc}\rho-1&0\\0&\rho+1\ea\right)\label{e-nf-k=01},
\mbox{
called {\bf S}$_1$-case},
\eqn
or the form,
\bqn
A=\left(\ba{cc}\eta&1\\0&\eta\ea\right),~~~
B=\left(\ba{cc}\rho+1&0\\0&\rho-1\ea\right)\label{e-nf-k=02},
\mbox{
called {\bf S}$_{-1}$-case}.
\eqn
\el
%%%%%%%%%%%%%%%%%%%%%%%
{\bf Proof of Lemma \ref{l-nf-kneq0} and Lemma \ref{l-nf-k=0}}.
We can always find a system of coordinates such that,
\bqn
&&A=\left(\ba{cc}\lam&1\\0&\lam\ea\right),~~~B=\left(\ba{cc}a&b\\c&d\ea\right),
~~a,b,c,d\in\R.
\eqn
In this case the discriminant of $B$ is $\delta=(a-d)^2+4 bc$. We have
\bqn
[A,B]=\left(\ba{cc}c&{d-a}\\0&-c\ea\right).
\eqn
Notice
moreover that $c=Tr(AB)-\frac12Tr(A)Tr(B)$. Hence, according with \r{e-k}, $c=0$ iff $k=0$.
\\\\
%%%%%%%%%%%%%%%%%%%%%%%%%%
{\bf Case $c\neq0$}. First notice that in this case $[A,B]\neq0$.
Consider the transformation
\bqn
&&T=\left(\ba{cc}1&\frac{a-d}{2c}\\0&1\ea\right).
\eqn
Then
\bqn
A':=T^{-1}AT=A,~~~~~~~~~
B':=T^{-1}BT=\left(\ba{cc}\frac{a+d}{2}&\frac{\delta}{4c}\\c&\frac{a+d}{2}
\ea\right).
\eqn
we have the following:
\bi
\i
If $\delta=0$, according to Definition \ref{d-invariants}, then $c=k$,
$\lam=\eta$, $\rho=(a+d)/2$, and with the change of notation
$B'\to B$ we get the normal forms \r{e-nA} and \r{e-nB}.
\i If $\delta\neq0$,  we define
\bqn
A'':=\para A'=\left(\ba{cc}\lam\para&\para\\0&\lam\para
\ea\right),~~~
B'':=\frac{2}{\sqrt{|\delta|}}B'=
\left(\ba{cc}\frac{a+d}{\sqrt{|\delta|}}&\frac{\sqrt{|\delta|}}{2c}\\\frac{2c}{\sqrt{|\delta|}}
&\frac{a+d}{\sqrt{|\delta|}}
\ea\right).
\eqn
Consider the transformation
\bqn
\llabel{T'}
&&T'=\left(\ba{cc}\frac{\sqrt{2}}{|\delta|^{1/4}}&0\\0&
\frac{|\delta|^{1/4}}{\sqrt{2}}\ea\right).
\eqn
We have
\bqn
A''':=(T')^{-1}A''T'=\left(\ba{cc}\lam\para&1\\0&\lam\para
\ea\right),~~~
B''':=(T')^{-1}B''T'=
\left(\ba{cc}\frac{a+d}{\sqrt{|\delta|}}&\frac{\delta}{4c}\\\frac{4c}{|\delta|}
&\frac{a+d}{\sqrt{|\delta|}}
\ea\right).
\eqn
According to Definition \ref{d-invariants},
and with the change of notations $A'''\to A$, $B'''\to B$ we get the
normal forms \r{e-nA} and
\r{e-nB}. Lemma \ref{l-nf-kneq0} is proved.
\ei
%%%%%%%%%%%%%%%%%%%%%%%%%%%%%%%%%%%%%%%%%%%%%%%%%%%%%%%%%%%%%%%%%%%%%%%%%%%%%%%%%
{\bf Case $c=0$}. In this case {\bf (H0)} implies $\delta=(a-d)^2>0$.
Consider the transformation:
\bqn
&&T=\left(\ba{cc}1&\frac{-b}{a-d}\\0&1\ea\right).
\eqn
We have:
\bqn
A':=T^{-1}AT=A,~~~~~~~~~
B':=T^{-1}BT=\left(\ba{cc} a&0\\0&d
\ea\right).
\eqn
We have two cases
\bi
\i $d>a$. In this case from $a+d=Tr(B)$ and $d-a=\sqrt{\delta}$, we get
$$a=\frac12(Tr(B)-\sqrt{\delta}),~~~d=\frac12(Tr(B)+\sqrt{\delta}).$$
Define
\bqn
A'':=\paraa A'=\left(\ba{cc}\lam\paraa&\paraa\\0&\lam\paraa
\ea\right),~~~
B'':=\frac{2}{\sqrt{\delta}}B=
\left(\ba{cc}    \rho-1&0\\0
&\rho+1
\ea\right).
\eqn
Define now $T'$ as in \r{T'}.
%\bqn
%\llabel{T'}
%&&T'=\left(\ba{cc}\frac{\sqrt{2}}{\delta^{1/4}}&0\\0&
%\frac{\delta^{1/4}}{\sqrt{2}}\ea\right).
%\eqn
Then
\bqn
A''':=(T')^{-1} A''T =\left(\ba{cc}\eta&1\\0&\eta
\ea\right),~~~
B''':=(T')^{-1} B''T=B''.
\eqn
With the change of notations $A'''\to A$, $B'''\to B$ we get the normal
forms
\r{e-nf-k=01}.
\i Similarly if  $d<a$, we get the normal forms  \r{e-nf-k=02}.
\ei
Lemma \ref{l-nf-k=0} is proved. \quadp\\\\
%%%%%%%%%%%%%%%%%%%%%%%%%%%%%%%%%%%%%%%%%%%%%%%%
The following Lemma can be directly checked.
\bl
\llabel{l-autovalori}
Under the condition {\bf (H0)}, we have the following:
\bi
\i $A$ has a unique eigenvalue $\eta<0$ and its corresponding eigenvector is $(1,0)$.
\i if $k\neq0$ and $\delta>0$ (i.e., in the {\bf R}$_{1}$-case),
 then the eigenvalues of $B$ are $\rho+1$ and $\rho-1$,
corresponding respectively to the eigenvectors $(1,k)$ and $(1,-k)$, and we have $\rho<-1$.
\i if $k\neq0$ and $\delta<0$, (i.e., in the   {\bf R}$_{-1}$-case) then
the
eigenvalues of $B$ are
$\rho+i$ and $\rho-i$, and we have $\rho<0$. In this case the integral curves of the vector field $Bx$ are elliptical
spirals rotating counter-clockwise if $k>0$ (clockwise if $k<0$).
\i if $k\neq0$ and $\delta=0$ (i.e., in the {\bf R}$_{0}$-case), then $B$
has a unique eigenvalue
$\rho<0$, and it corresponds
to the eigenvector $(0,1)$.
\i if $k=0$, (i.e., in the  {\bf S}-case) then the eigenvalues of $B$ are
$\rho+1$ and
$\rho-1$, with $\rho<-1$.
\bi
\i In the case
of the normal form \r{e-nf-k=01} (i.e., in the {\bf S}$_{1}$-case)
they correspond respectively
to the eigenvectors $(0,1)$ and $(1,0)$.
\i In the case of the normal form \r{e-nf-k=02}
(i.e., in the {\bf S}$_{-1}$-case), they
correspond respectively to the eigenvectors
$(1,0)$ and $(0,1)$.
\ei
\ei
\el

\section{Main Results}
\llabel{s-main}
In this section we state our stability conditions. First we need to define
some functions of the invariants $\eta,$ $\rho,$ $k$ defined  in
Definition \ref{d-invariants}.
Set
\bqn
\Delta&=& k^2 - 4 \eta \rho k + sign(\delta) 4\eta^2.
\label{DELTA}
\nn
\eqn
By direct computation one gets that if $k=2\eta\rho$ then
$\Delta=-4 \det A \det B<0$. It follows
\bl
\label{l-ovvio}
Assume {\bf (H0)}. Then $\Delta\geq0$ implies $k\neq2\eta\rho$.
\el
%%%%%%%%%%%%%%%%%%%%%%%%%%%%%%%%%%%%%%%%%%%%%%%%%%%%%%%%
Moreover, when $\Delta>0$ and $k<0$, define
\bqn
{\cal R}= \left\{
%%%
\ba{l} |(\frac{-k+\sqrt{\Delta}}{-k-\sqrt{\Delta}})\frac{2 k \rho^2
- sign(\delta)(Tr(AB)+ \sqrt{\Delta})}{2 k \sqrt {\det B}
 (\rho-sign(\delta)\frac{\eta}{k})}|\exp(\frac{\sqrt{\Delta}}{k}+\rho
  \theta_{sign(\delta)} ), \ \ {\rm if}
 \ \ \rho-sign(\delta)\frac{\eta}{k}\neq 0, \\
%%%
\frac{-2 \eta}{\sqrt{k^2+
\eta^2}}\exp(\frac{\sqrt{\Delta}}{k}+\rho \theta_{-1} ), \ \ {\rm
if} \ \ \rho-sign(\delta)\frac{\eta}{k} = 0 \ \  ({\rm which \ \
implies} \ \ sign(\delta)=-1).
\ea\right.
\eqn
where
\bqn
\theta_{-1} &=& \left\{\ba{ll}\arctan \frac{\sqrt{\Delta}}{k(\rho +
\frac{\eta}{k})+
\eta}
&~~\mbox{ if }~~k(\rho +
\frac{\eta}{k})+\eta\neq0\\
\pi/2&~~\mbox{ if }~~k(\rho +
\frac{\eta}{k})+\eta=0,
\ea\right. \nn\\
\theta_1 &=& \mbox{arctanh}\frac{\sqrt{\Delta}}{k(\rho - \frac{\eta}{k})-
\eta},\nn\\
\theta_0 &=&\frac{\sqrt{\Delta}}{k \rho}. \nn
\eqn
Notice that when $k<0,$ then $k(\rho - \frac{\eta}{k})- \eta >0 \ \
\text{and} \ \ k \rho >0$. Hence $\theta_1$ and $\theta_0$ are well
defined. \\\\
%%%%%%%%%%%%%%%%
The following Theorem  states the main result of the paper. The letters
{\bf A.},
{\bf B.}, and
{\bf C.} refer again to the cases described in the introduction and in Figure \ref{f-stabbb}.
Recall Lemma \ref{l-ovvio}.
%%%%%%%%%%%%%%%%%%%%%%%%%%%%%%%%%%%%%%%%%%%%%%%%%%%%%%%%%%%
\begin{Theorem}
\label{t-main}
Assume {\bf (H0)}. We have the following stability conditions for the
system \r{1}.
\bd
    \item[A.] If $\Delta <0,$ then the system is GUAS.

    \item[B.] If $\Delta >0,$ then:
    \bd
        \item[B1.] if $k> 2\eta \rho$, then the system is unbounded,
        \item[B2.]if $k< 2\eta \rho$, then
\bi
\item in the regular case ($k \neq 0$), the system is GUAS,
uniformly stable (but not GUAS) or unbounded respectively if
                $$ {\cal R} <1,{\cal R} =1, {\cal R} >1 .$$
\item In the singular case ($k=0$), the system is GUAS.
\ei
\ed
\item[C.] If $\Delta = 0,$ then:
    \bd
    \item[C1] If  $k>2\eta \rho,$ then the system is uniformly stable (but not
    GUAS),

    \item[C2] if $k<2\eta \rho,$ then the system is GUAS.
    \ed
      \ed
\end{Theorem}
%%%%%%%%%%%%%%%%%%%%%%%%%%%%%%%%%%%%%%%%%%%%%%%%%%%%%%%%%%%%%%%%%%%%%%%%
%%%%%%%%%%%%%%%%%%%%%%%%%%%%%%%%%%%%%%%%%%%%%%%%%%%%%%%%%%%%%%%%%%%%%%%%
\section{The set where the two vector fields are parallel}
\llabel{s-parallel}
As explained in the introduction, the stability conditions are obtained by
studying the locus in which the
two vector fields $Ax$ and $Bx$ are linearly dependent. This set is the
set $\qq$ of zeros of the quadratic form $Q(x):=Det(Ax,Bx)$.
Set $x=(x_1,x_2)$. We have the following:
\bqn
Q=\left\{\ba{ll}
\left(\rho-sign(\delta)\frac{\eta}{k}\right)x_2^2+k x_1x_2+(\eta
k)x_1^2&\mbox{ in the {\bf R}-case}\\
x_2\left((\rho+1)x_2+2\eta x_1\right)&\mbox{ in the {\bf S}$_{1}$ case }\\
x_2\left((\rho-1)x_2-2\eta x_1\right)&\mbox{ in the {\bf S}$_{-1}$ case. }
\ea\right.
\eqn
The discriminant of this quadratic form is the quantity $\Delta$ defined
in \r{DELTA}.
%%%%%%%%%%%%%%%%%%%%%%%%%%%%%%%%%%%%%%%%%%%%%%%%%%%%%%%%%%%%%%%%%%
Hence we have the following:
\bi
\i if $\Delta<0$ then $\qq=\{0\}$,
\i if $\Delta>0$ then $\qq$ is a pair of
(noncoinciding) straight lines passing through the origin,
\i if $\Delta=0$ then  $\qq$ is a single straight line passing through the
origin (the two straight lines of the previous case coincide).
\ei
%%%%%%%%%%%%%%%%%%%%%%%%%%%%%%%%%%%%%%%%%%%%%%%%%%%%%%%%%%%%%%%%
\brem
Notice that, under the condition {\bf (H0)}, in the regular case, the sign
of
$\Delta$
can be positive or negative, while in the singular case, since $k=0$ and
$\delta>0$ (cf. Lemma \ref{l-nf-k=0}), we have $\Delta=4\eta^2>0$.
In the singular case $\qq$ is the union of a  pair of
straight lines one of them coinciding with the $x_1$ axis (as it is
clear from the expression of $Q$).
\erem
%%%%%%%%%%%%%%%%%%%%%%%%%%%%%%%%%%%%%%%%%%%%%%%%%%%%%%%%%%%%%%%%%%
In the following, under the assumption $\Delta\geq0$,  we give the
explicit expression of the angular coefficients of the two
(possibly coinciding) straight lines whose union is $\qq$.
\bd
\i[{\bf R}$_{-1}$-case]  (i.e., the case in which $B$ has non real
eigenvalues)
\bqn
\mbox{ if }\rho+\eta/k\neq0\mbox{ then }
&\displaystyle
m^\pm=\frac{-k\pm\sqrt{\Delta}}{2(\rho+\eta/k)}=\frac{-k\pm\sqrt{k^2-4\eta
k(\rho+\eta/k)}}{2(\rho+\eta/k)}\label{m+-1},\\
\mbox{ if }\rho+\eta/k=0\mbox{ then } &m^+=\infty,~~~
m^-=-\eta.\label{m--1}
\eqn
\i[{\bf R}$_{1}$-case] (i.e., the regular case in which $B$ is
diagonalizable and has real eigenvalues)
\bqn
\mbox{ if }\rho-\eta/k\neq0\mbox{ then }
&\displaystyle
m^\pm=\frac{-k\pm\sqrt{\Delta}}{2(\rho-\eta/k)}=\frac{-k\pm\sqrt{k^2-4\eta
k(\rho-\eta/k)}}{2(\rho-\eta/k)}\label{m++1},\\
\mbox{ if }\rho-\eta/k=0\mbox{ then } &m^-=\infty,~~~ m^+=-\eta.\label{m-+1}
\eqn
\i[{\bf R}$_{0}$-case] (i.e., the regular case in which $B$ is
nondiagonalizable)
\bqn
m^\pm=\frac{-k\pm\sqrt{\Delta}}{2\rho}=\frac{-k\pm\sqrt{k^2-4\eta
k\rho}}{2\rho}.
\eqnl{m00}

\i [{\bf S}-case] (i.e., the singular case)
\bqn
m^+=0,~~~m^-=-\frac{2\eta}{\rho+1}<0\mbox{ in the {\bf S$_{1}$}-case},\\
m^+=0,~~~m^-=\frac{2\eta}{\rho-1}>0\mbox{ in the {\bf S$_{-1}$}-case}.
\eqn
\ed
%%%%%%%%%%%%%%%%%%%%%%%%%%%%%%%%%%%%%%%%%%%%%%%%%%%%%%%
The following Proposition says that  if $\Delta\geq0$ and  the two  vector
fields have the same versus on a point of $\qq\setminus\{0\}$,
then this is the case all along $\qq\setminus\{0\}$.
More precisely it says that if $k<2\eta\rho$ (resp. $k>2\eta\rho$) then they
have the same (resp. opposite) versus.
Recall Lemma \ref{l-ovvio}.
%%%%%%%%%%%%%%%%%%%%%%%%%%%%%%%%%%%%%%%%%%%%%%%%%%%%%%%%%%%%%%%%%%%%%

\bp
\label{p-or} Assume {\bf (H0)}, $\Delta\geq0$ , and let
$\qq= D^{+}\cup D^{-}  \ \ \text{where} \ \
D^{\pm}=\{(h,m^{\pm}h)\in\mathbb{R}^2, \ h\in \mathbb{R}\}$. Let
us define $\al^\pm$ by $Bx = \al^\pm Ax$ for $x\in D^\pm
\setminus\{0\}$. Then
  \bi
    \i $\al^+ \al^-  =  \frac{\det B}{\det A}>0$,
    \i $\al^+ + \al^- = \frac{2\eta\rho -k}{\det A}.$
  \ei
In other words $sign(\al^+)=sign(\al^-)=sign(2\eta\rho-k)$.
\ep
%%%%%%%%%%%%%%%%%%%%%%%%%%%%%%%%%%%%%%%%%%%%%%%%%%%%%%%%%%%%%%%%%%
\proof
Let us start with the regular case and $\Delta>0$.   In this case it is easy to check that
$\al^{\pm}=
\frac{k+\rho m^{\pm}}{\eta m^{\pm}}$. Hence $\al^+\al^- = \frac{k^2+k\rho(m^+ +
m^-)+\rho m^+m^-}{\eta^2m^+m^-}$. Using the fact that $m^+ + m^-
= \frac{-k}{\chi}$ and $m^+m^- = \frac{\eta
k}{\chi}$, where $\chi:=\rho-sign(\delta) \frac{\eta}{k}
$, we get
$$\al^+\al^- = \frac{k(\chi - \rho)+\rho^2 \eta}{\eta^3} = \frac{\rho^2
\eta - sign(\delta)\eta}{\eta^3} = \frac{\rho^2-
sign(\delta)}{\eta^2}= \frac{\det B}{\det A}.$$\\
Similarly we have $$\al^+ + \al^- = \frac{k(m^++m^-)+ 2\rho
m^+m^-}{\eta m^+ m^-}= \frac{-k^2+ 2\rho \eta k}{\eta^2
k}=\frac{2\eta\rho -k}{\det A}.$$
%%%%%%%%%%%%%%%%%%%%%%%%%%%%%%%%%%%
In the regular case, with $\Delta=0$, we have $D^+=D^-$, $m^+=m^-=:m^0$ and
$\alpha^+=\alpha^-=:\alpha^0.$ Now $Bx_1^{0} =
\al^{0}Ax_1^{0}$ implies that
$$
\alpha^0= \frac{k+\rho m^0}{\eta
m^0}= \frac{k+\rho \frac{-k}{2 \chi}}{\eta
 \frac{-k}{2 \chi}} =\frac{\rho-2 \chi }{\eta}.
$$
Since $\Delta= k^2-4\eta k \chi=0$ we have $2 \chi =
 \frac{k}{2\eta}$ and  $\alpha^0=
\frac{\rho - \frac{k}{2\eta}}{\eta} = \frac{2\eta\rho
-k}{2\eta^2}.$ Thus
$$
2\alpha^0 = \frac{2\eta\rho -k}{\det A} \
\  \text{and} \ \ \ (\alpha^0)^2
=\frac{\rho^2-sign(\delta)}{\eta^2}=\frac{\det B}{\det A}
$$
%%%%%%%%%%%%%%%%%%%%%%%%%%%%%%%%%%%%%%%%%%
In the singular case $k=0$,  we have
$\Delta>0$ and
\bqn
&&\qq =\{(x_1,x_2)\mid x_2=0\}\cup\{(x_1,x_2)\mid x_2=
\frac{-2\eta}{\rho +1}x_1\}
\text{~~~or}\nn\\
&&\qq=\{(x_1,x_2)\mid x_2=0\}\cup\{(x_1,x_2)\mid x_2=
\frac{2\eta}{\rho
-1}x_1\}.\nn
\eqn
An easy computation show that
$$\alpha^+ \alpha^-=\frac{(\rho-1)(\rho+1)}{\eta^2}=\frac{\det B}{\det A}
\ \ \text{and} \ \ \alpha^+ +
\alpha^-= \frac{2\rho}{\eta}= \frac{2 \eta\rho}{\eta^2}=
\frac{2\eta\rho}{\det A}.$$\quadp
%%%%%%%%%%%%%%%%%%%%%%%%%%%%%%%%%%%%%%%%%%%%%%%%%%%%%%%%%%%%
\bdeff
If $\Delta\geq0$ and $k<2\eta\rho$
(resp.  $k<2\eta\rho$) we say that $\qq$ is ``direct''
(resp. ``inverse'').
\edeff
The following Proposition says, roughly speaking, that on $\qq\setminus
\{0\}$ the
vector
field $Ax$
points ``clockwise''
(see Figure \ref{f-rotation}).
\bp
\llabel{p-rotation}
Assume  {\bf (H0)} and $\Delta\geq0$.
Let $\partial_\theta$ be the vector field defined by $\partial_\theta(x_1,x_2)=(-x_2,x_1)$.
Then in the normal forms of Section \ref{s-normal}, we have
$Ax\cdot\partial_\theta\leq0$ for every $x\in
\qq\setminus\{0\}$.
\ep
\proof
Let $x^+\in D^+$, where $D^+$ is defined in Proposition \ref{p-or}. Then
$x^+=(h,m^+h)$ for some $h\in \R\setminus\{0\}$.
It
follows $Ax^+\cdot\partial_\theta(x^+)=-(m^+)^2h^2$.
Similarly if
$x^-=(\bar h,m^-\bar h)$ for some $\bar h\in \R\setminus\{0\}$,
we have $Ax^-\cdot\partial_\theta(x^-)=-(m^-)^2\bar h^2$.\quadp
\ppotR{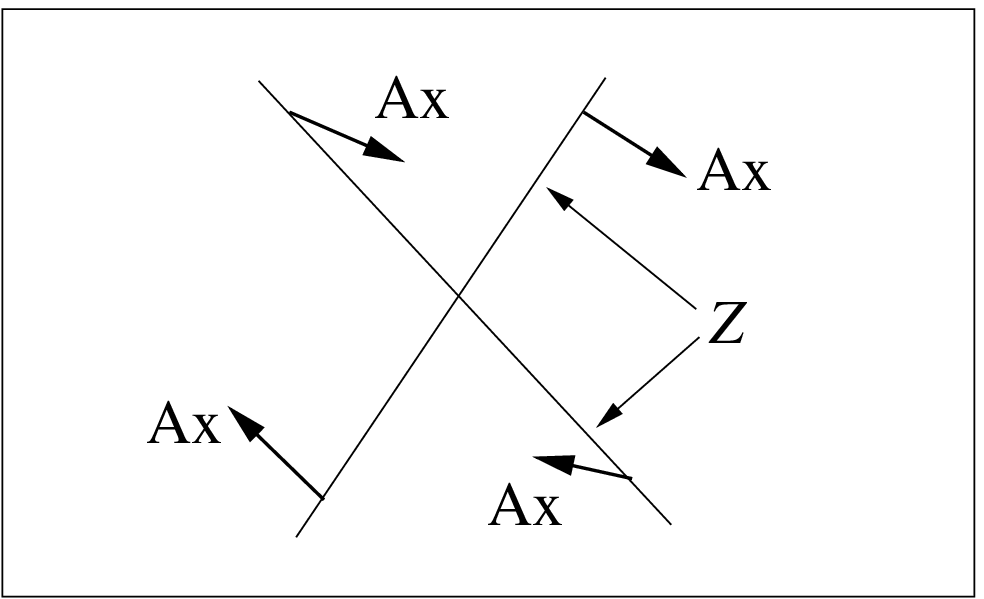}{}{7}

%%%%%%%%%%%%%%%%%%%%%%%%%%%%%%%%%%%%%%%%%%%%%%%%%%%%%%%%%%%%%%
%%%%%%%%%%%%%%%%%%%%%%%%%%%%%%%%%%%%%%%%%%%%%%%%%%%%%%%%%%%%%%
%%%%%%%%%%%%%%%%%%%%%%%%%%%%%%%%%%%%%%%%%%%%%%%%%%%%%%%%%%%%%%
\section{General Stability Conditions}
\llabel{s-general}
%%%%%%%%%%%%%%%%%%%%%%%%%%%%%%%%%%%%%%%%%%%%%%%%%%%%%%%%%%%%%%%
In this Section we study some basic stability conditions.

The following proposition (from which it follows {\bf A.} of Theorem \ref{t-main}) expresses the
idea that if the two
vector fields are linearly dependent just at the origin, then the system is always
GUAS. For the proof one can follow the same steps used in the case of two
diagonalizable matrices \cite{sw-1}. It is also a particular case of the
same result proved for nonlinear systems in \cite{sw-nonlinear}.
%%%%%%%%%%%%%%%%%%%%%%%%%%%%%%%%%
\bp
\label{p-delta<0}
Assume {\bf H0}. If $\Delta<0$ then the system is GUAS.
\ep
%%%%%%%%%%%%%%%%%%%%%%%%%%%%%%%%%%%%%%%%%%%%%%%%%%%%%%%%%%%%%%%%%
\subsection{The inverse case}
\label{ss-inverse}
The next proposition (from which it follows {\bf B1.} of Theorem \ref{t-main}) says that
if the vector fields are parallel along two (noncoinciding) straight lines
and point in opposite directions, then
there are trajectories going to infinity.
\bp
\label{p-2etarho}
Assume {\bf (H0)}. If  $\Delta>0$ and $k>2\eta\rho$ then
the system is unbounded. In the degenerate case in which $\Delta=0$ and
$k>2\eta\rho$ then the system is uniformly stable but not GUAS.
\ep
%%%%%%%%%%%%%%%%%%%%%%%%%%%%%%%%%%%%%%%%%%%%%%%%%%%%%%%%%%%%%%%%%
\proof  To prove the first statement, recall that we can work with the
convexified system. Consider the following
constant control $$u_0 = \frac{\al^+ \al^- - \frac{1}{2}(\al^+ +
\al^-)}{1+\al^+ \al^- -(\al^+ + \al^-) }.$$ From Proposition
\ref{p-or} we have $\al^+ + \al^- <0 \  \ \text{and}  \  \
\al^+ \al^- >0,$ which implies that $u_0 \in ]0,1[.$ A simple computation show that: $$\det
M(u_0)= \det(u_0 A +(1-u_0)B) = \frac{-\Delta}{4 \det A(1+\al^+ \al^- -(\al^+ +
\al^-)} <0, $$ since $\Delta>0.$ Hence $M(u_0)$ has a positive real
eigenvalue and the system is unbounded.\\ \\
%%%%%%%%%%%%%%%%%%%%%%%%%%%%%%%%%%%%%%%%%%%%%%%%%
To prove the second statement, observe that if $\Delta=0$ and $k> 2 \eta \rho$, then   $\qq=
D^{+}=D^- = \{(x_1,x_2) \ | \ x_2 + 2\eta x_1 =0 \},$ and
thanks to Proposition \ref{p-or}, $\al^+ = \al^-=\al^0 <0$.  Consider
again the matrix $$M(u_0)=u_0 A + (1-u_0)B \ \ \text{where} \ \
u_0 = \frac{(\al^0)^2 -\al^0}{1- 2\al^0 + \al^0{^2}} =
\frac{\al^0}{\al^0-1}.$$ Then
$Tr(M(u_0))= \frac{2(\al^0 \eta - \rho)}{\al^0-1}<0$ and
$\det M(u_0)=0.$ Therefore  for $u_0$ the system is not asymptotically
stable. \\
To prove that the system is uniformly stable we show that
$Ax$ and $Bx$ admit the following  common (non strict) Lyapunov function,
$$V(x)= V(x_1,x_2)=x_1^2 + \frac{x_2^2}{4 \eta^2}.$$
Let us prove that
\bqn
\nabla V(x)Ax &=& \frac{1}{2 \eta}(x_2 + 2\eta
x_1)^2\leq 0 \label{lu1}\\
\nabla V(x)Bx & =& \frac{\rho}{2 \eta^2}(x_2 + 2\eta x_1)^2\leq 0.
\label{lu2}
\eqn
To prove \r{lu1} observe that $\nabla
V(x)Ax= 2\eta x_1^2 + 2 x_1x_2 + \frac{1}{2\eta}x_2^2 = \frac{1}{2
\eta}(x_2 + 2\eta x_1)^2.$
To prove \r{lu2}, observe that $$\nabla V(x)Bx= 2\rho x_1^2 + (\frac{2
sign(\delta)}{k}+ \frac{k}{2 \eta^2})x_1 x_2+
\frac{\rho}{2\eta^2}x_2^2=2\rho x_1^2 + (\frac{4
sign(\delta)\eta^2+k^2}{2\eta^2 k})x_1 x_2+
\frac{\rho}{2\eta^2}x_2^2.$$ But $\Delta = k^2 - 4\eta k \rho + 4
\eta^2 sign(\delta)=0$ implies that  $k^2  + 4 \eta^2
sign(\delta)= 4\eta k \rho$. Then  $$\nabla V(x)Bx= 2\rho x_1^2 +
\frac{4 \eta k \rho}{2 \eta^2 k }x_1 x_2+
\frac{\rho}{2\eta^2}x_2^2 = 2\rho x_1^2 + \frac{2 \rho}{ \eta}x_1
x_2+ \frac{\rho}{2\eta^2}x_2^2 = \frac{\rho}{2 \eta^2}(x_2 + 2\eta
x_1)^2 .$$ It follows that  $V(x)$ is a common Lyapunov function and the system
is uniformly stable.
\quadp\\\\
%%%%%%%%%%%%%%%%%%%%%%%%%%%%%%%%%%%%%%%%%%%%%%%%%%%%%%%%%%%%%%%%%%%%%%%%%%%
Notice that in the {\bf S}-case, Proposition \ref{p-2etarho} never
apply because $\qq$ is always ``direct''.
%%%%%%%%%%%%%%%%%%%%%%%%%%%%%%%%%%%%%%%%%%%%%%%%%%%%%%%%%%%%%%%%%%%
%%%%%%%%%%%%%%%%%%%%%%%%%%%%%%%%%%%%%%%%%%%%%%%%%%%%%%%%%%%%%%%%%%%
%%%%%%%%%%%%%%%%%%%%%%%%%%%%%%%%%%%%%%%%%%%%%%%%%%%%%%%%%%%%%%%%%%%
\subsection{The direct case: the worst-trajectory}
In the case in which  $\Delta\geq0$ and  $k<2\eta\rho$, the
stability of \r{1} can be reduced to the study of a single trajectory
called the ``worst-trajectory''.
%%%%%%%%%%%%%%%%%%%%%%%%%%%%%%%%%%%%%%%%%%%%%%%%%%%%%%%%%%%%%%%%%%%%
\bdeff
Assume   {\bf (H0)}, $\Delta\geq0$ and  $k<2\eta\rho$.
Fix $x_0\in\R^2\setminus\{0\}$. The worst-trajectory
$\g_{x_0}$ is the trajectory of \r{1}, based at
$x_0$, and  having the following
property. At each time $t$, $\dot \g_{x_0}(t)$ forms the
smallest angle (in absolute value) with the
(exiting) radial direction (see Figure \ref{f-worst}).
\llabel{d-worst}
\edeff
As explained in the introduction, (see also \cite{sw-1,sw-lyapunov}) the
worst-trajectory
switches among the two vector fields on the set $\qq$.
%%%%%%%%%%%%%%%%%%%%%%%%%%%%%%%%%%%%%%%%%%%%%%%%%%%%%%%%%%%%%%%%%%%%%
The following Lemma (whose proof is a consequence of the arguments used in \cite{sw-1}) reduce the
problem of stability of \r{1} to the problem of
the stability of the worst-trajectory.
\bl
\llabel{l-w}
Assume {\bf (H0)}, $\Delta>0$, and $k<2\eta\rho$.
Fix $x_0\in\R^2\setminus\{0\}$.
Then the system \r{1} is GUAS (resp. uniformly stable but not GUAS, resp. unbounded)
if and only if $\g_{x_0}$ tends to the origin (resp. is periodic, resp. tends to infinity).
Moreover when $\g_{x_0}$  is periodic or tends to infinity then it rotates around the origin switching an
infinity number of times. In the degenerate case in which  $\Delta=0$, and $k<2\eta\rho$, then the system
is GUAS.
\el
%%%%%%%%%%%%%%%%%%%%%%%%%%%%%%%%%%%%%%%%%%%%%%%%%%%%%%%%%%%%%%%%%
\brem
The last statement can be be proved either by using worst-trajectory
type arguments or by building a common
quadratic Lyapunov function and it implies statement {\bf C2.} of Theorem
\ref{t-main}.
\erem
%%%%%%%%%%%%%%%%%%%%%%%%%%%%%%%%%%%%%%%%%%%%%%%%%%%%%%%%%%%%%%%%%%%%%%%%%%
\brem
\llabel{r-rotation}
The explicit construction of the worst trajectory in the case $\Delta>0$ and $k<2\eta\rho$
(i.e., the proof of case {\bf B2.} of Theorem \ref{t-main})
is done in
the next section.
The worst-trajectory can rotate or not around the origin.
As it will be  clear next, if  $B$ has complex eigenvalues then it always rotates. One of the statement
of Lemma \ref{l-w} is that if the worst-trajectory does
not rotate then the system is always GUAS. Notice that, if the worst-trajectory
rotates around the
origin, then it rotates clockwise, thanks to Proposition
\ref{p-rotation}.
\erem
%%%%%%%%%%%%%%%%%%%%%%%%%%%%%%%%%%%%%%%%%%%%%%%%%%%%%%%%%%%%%%%%%%%%%%%%%%%%%%%%%%%%%%%%
%%%%%%%%%%%%%%%%%%%%%%%%%%%%%%%%%%%%%%%%%%%%%%%%%%%%%%%%%%%%%%%%%%%%%%%%%%%%%%%%%%%%%%%%
%%%%%%%%%%%%%%%%%%%%%%%%%%%%%%%%%%%%%%%%%%%%%%%%%%%%%%%%%%%%%%%%%%%%%%%%%%%%%%%%%%%%%%%%
\section{Construction of the worst-trajectory}
\llabel{s-worst}
In this section we prove {\bf B2.} of Theorem \ref{t-main}, separately for the cases
{\bf R}$_{-1}$,
{\bf R}$_{1}$,
{\bf R}$_{0}$, and
{\bf S}.
%%%%%%%%%%%%%%%%%%%%%%%%%%%%%%%%%%%%%%%%%%%%%%%%%%%%%%%%%%%%%%%%%%%%%%%%%%%%
\subsection{The R$_{-1}$-case}
In the {\bf R}$_{-1}$-case we have $\delta<0$, therefore   the
matrix $B$ has non real eigenvalues.

Notice that under the conditions $\delta<0$ and $0<k<2\eta\rho$ we have
$\Delta=k^2-4\eta k(\rho+\eta/k)=k(k-2\eta\rho)-2k\eta\rho-4
\eta^2<0$. Therefore, since we are interested to the case in which
$\Delta>0$, we may assume $k<0$. Recall formulas \r{m+-1} and \r{m--1}. We
have the following
cases depicted in Figure \ref{f-r-1}, for $x_1\geq0$.
\ppotR{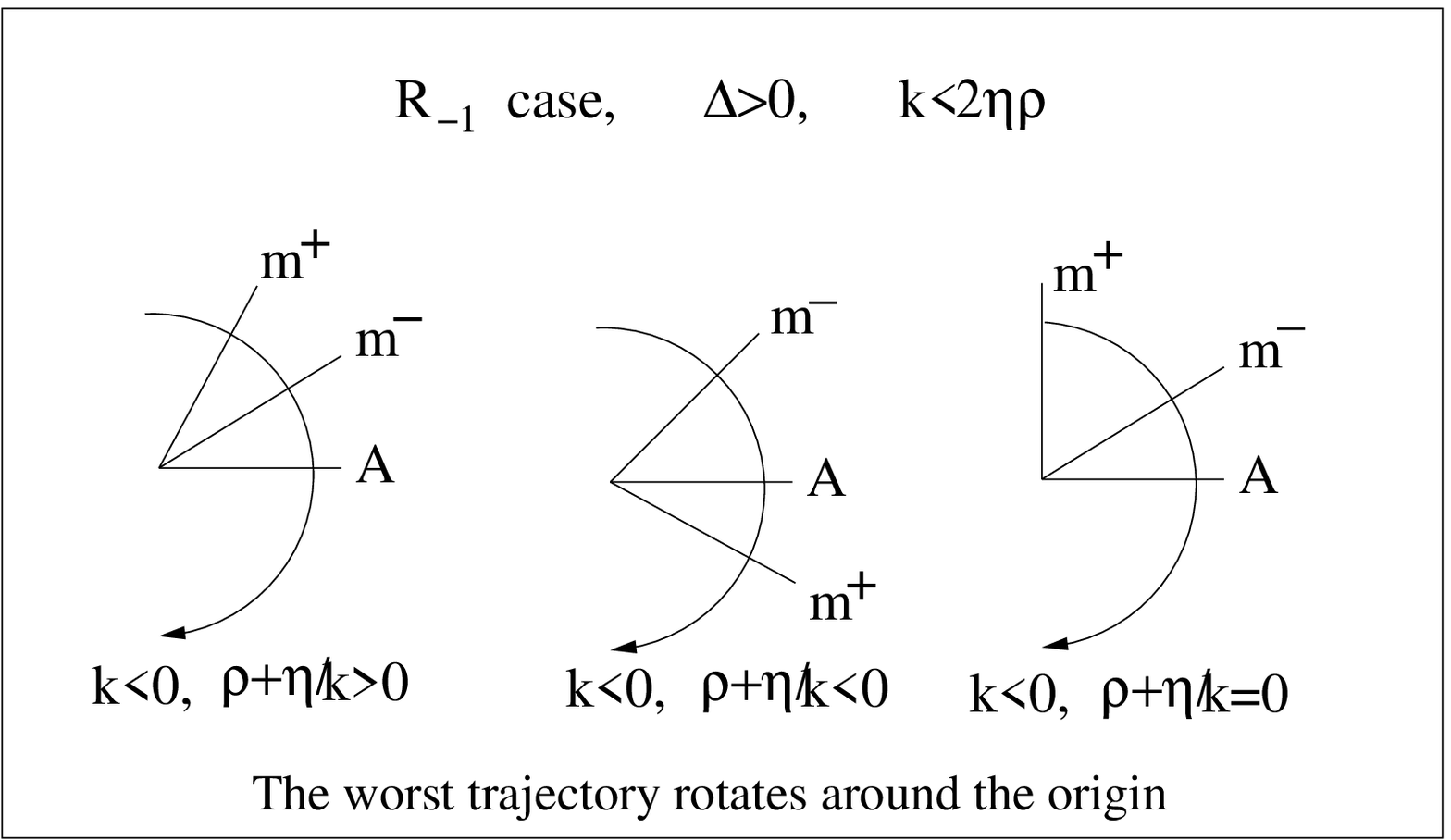}{}{11}
\bi
\i $k<0$
and $\rho+\eta/k>0$. In this case $\sqrt\Delta<|k|$. It follows
that $m^+>m^->0.$
\i $k<0$ and $\rho+\eta/k<0$. In this case
$\sqrt\Delta>|k|$. It follows that $m^->0>m^+$.
\i $k<0$ and $\rho+\eta/k=0$.
In this case $m^+=\infty$, $m^-=-\eta>0$.
\ei
Fix
$x_0=(1,m^+)\in D^+$. In all cases, the worst-trajectory
$\gamma_{x_0}$ rotates clockwise (cf. Remark \ref{r-rotation}) and
is a concatenation of integral curves of $Ax$ (from the line $D^+$
of equation $x_2= m^+ x_1$ to line $D^-$ of equation $x_2 = m^-
x_1,$) and integral curves of $Bx$ otherwise. See Figure
\ref{f-ruota-R-1}.
%%%%%%%%%%%%%%%%%%%%%%%%%%%%%%%%%%%%%%%%
\ppotR{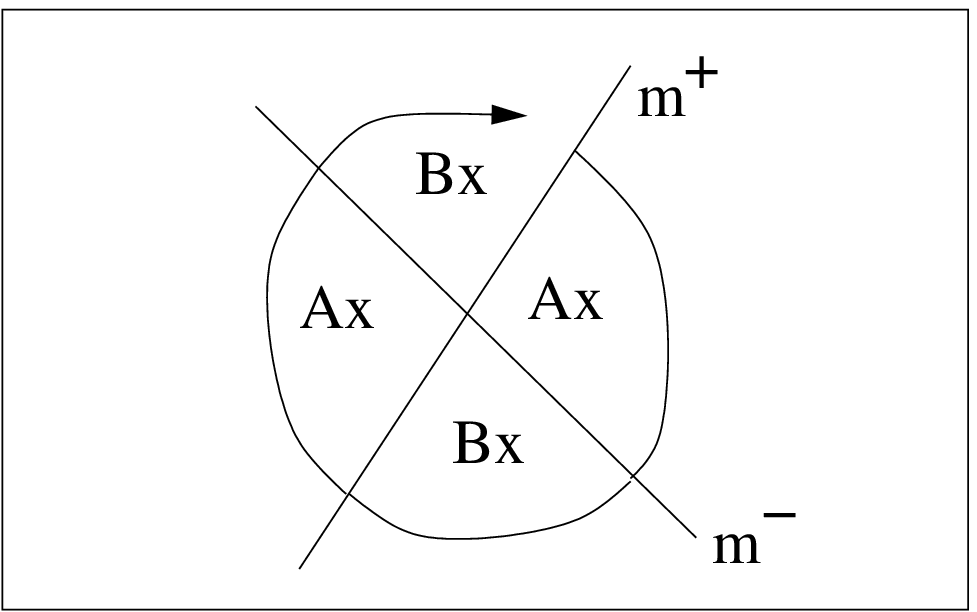}{}{6}
%%%%%%%%%%%%%%%%%%%%%%%%%%%%%%%%%%%%%%%%
Let $t_1$ be the time
needed to  an integral curve of $Ax$ to steer a point from $D^+$ to
$D^-$. Since
\bqn
\exp(At)=\exp(\eta t)\left (
\begin{array}{cc}
1 & t \\
0& 1
\end{array}
\right ),
\eqnl{At}
we have $$t_1= \frac{m^+-m^-}{m^+m^-}=\frac{\sqrt{\Delta}}{\eta
k}.$$ Let $t_2$ be the time needed to  an integral curve of $Bx$
to steer a point from $D^-$ to $D^+$. Since
$$
\exp(Bt) = \exp(\rho t) \left (
\begin{array}{cc}
\cos t & \frac{-\sin t}{k}\\
k \sin t & \cos t
\end{array}
\right ),$$
we have
$$
t_2 = \left\{\ba{ll}
\arctan \frac{\sqrt{\Delta}}{k(\rho + \frac{\eta}{k})+ \eta}&\text{ if } \ \ k(\rho +
\frac{\eta}{k})+ \eta \neq 0\\
\frac{\pi}{2} &\text{ if } \ \ k(\rho +
\frac{\eta}{k})+ \eta = 0.\ea\right.
$$
Notice that in the case $\rho +
\frac{\eta}{k}= 0,$ we have the following simple formulas for the
switching time: $t_1 = \frac{-1}{\eta}$ and $t_2 = \arctan\frac{-k}{\eta}$.
Fix $x_0\in D^+$ and let ${\cal R}$ be the ratio between the norm of
$\g_{x_0}$ after half turn (i.e., after time $t_1+t_2$)
and the norm of $\g_{x_0}(0)$, i.e.,
\bqn {\cal R}=\frac{\|\gamma_{x_0}(t_1+t_2)\|}{\|\gamma_{x_0}(0)\|}
\eqnl{R}
We have the following cases
\begin{itemize}
    \item $ \rho + \frac{\eta}{k} \neq 0.$ In this case
    $${\cal R}= \exp(\frac{\sqrt{\Delta}}{k} )\exp(\rho t_2 )
    |(1+\frac{m^+\sqrt{\Delta}}{k \eta})(\cos t_2 - \frac{m^-}{k}\sin t_2)|,$$ after simplification we get $${\cal R}
    =|(\frac{-k+\sqrt{\Delta}}{-k-\sqrt{\Delta}})\frac{2 k \rho^2 +
Tr(AB)+ \sqrt{\Delta}}{2 k \sqrt {\det B}
 (\rho+\frac{\eta}{k})}|\exp(\frac{\sqrt{\Delta}}{k} )\exp(\rho t_2),$$
i.e., $${\cal R}= |
    (\frac{-k+\sqrt{\Delta}}{-k-\sqrt{\Delta}})\frac{2 k \rho^2 -
    sign(\delta) (Tr(AB)+ \sqrt{\Delta})}{2 k \sqrt {\det B}
 (\rho-sign(\delta)\frac{\eta}{k})}|\exp(\frac{\sqrt{\Delta}}{k} + \rho \theta_{-1}).$$

\item $ \rho + \frac{\eta}{k}= 0$. In this case
 $${\cal R}= |\exp(\frac{\sqrt{\Delta}}{k} )
    \exp(\rho t_2 )(\cos t_2 + \frac{m^-}{k}\sin t_2)|=\frac{-2 \eta}{\sqrt{k^2+
\eta^2}}\exp(\frac{\sqrt{\Delta}}{k}+\rho \theta_{-1}).$$
\end{itemize}
The system is GUAS, uniformly stable (but not GUAS) or unbounded respectively if
${\cal R}<1$, ${\cal R}=1$, ${\cal R}>1$.
%%%%%%%%%%%%%%%%%%%%%%%%%%%%%%%%%%%%%%%%%
\subsection{The R$_{1}$-case}

In this section we study the R$_{1}$-case i.e., the case in which the
matrix $B$ has real eigenvalues. In this case since $\delta>0$, we have
$\Delta=k^2-4\eta k(\rho -\eta/k)$.
Again we restrict to the case $\Delta>0$, and $k<2\eta\rho$.
Set
\bqn
\chi:=\rho - \frac{\eta}{k}.
\eqn
Recall formulas \r{m++1} and \r{m-+1}. We have the following cases
depicted in Figure
\ref{f-r+1}, for $x_1>0$.
\ppotR{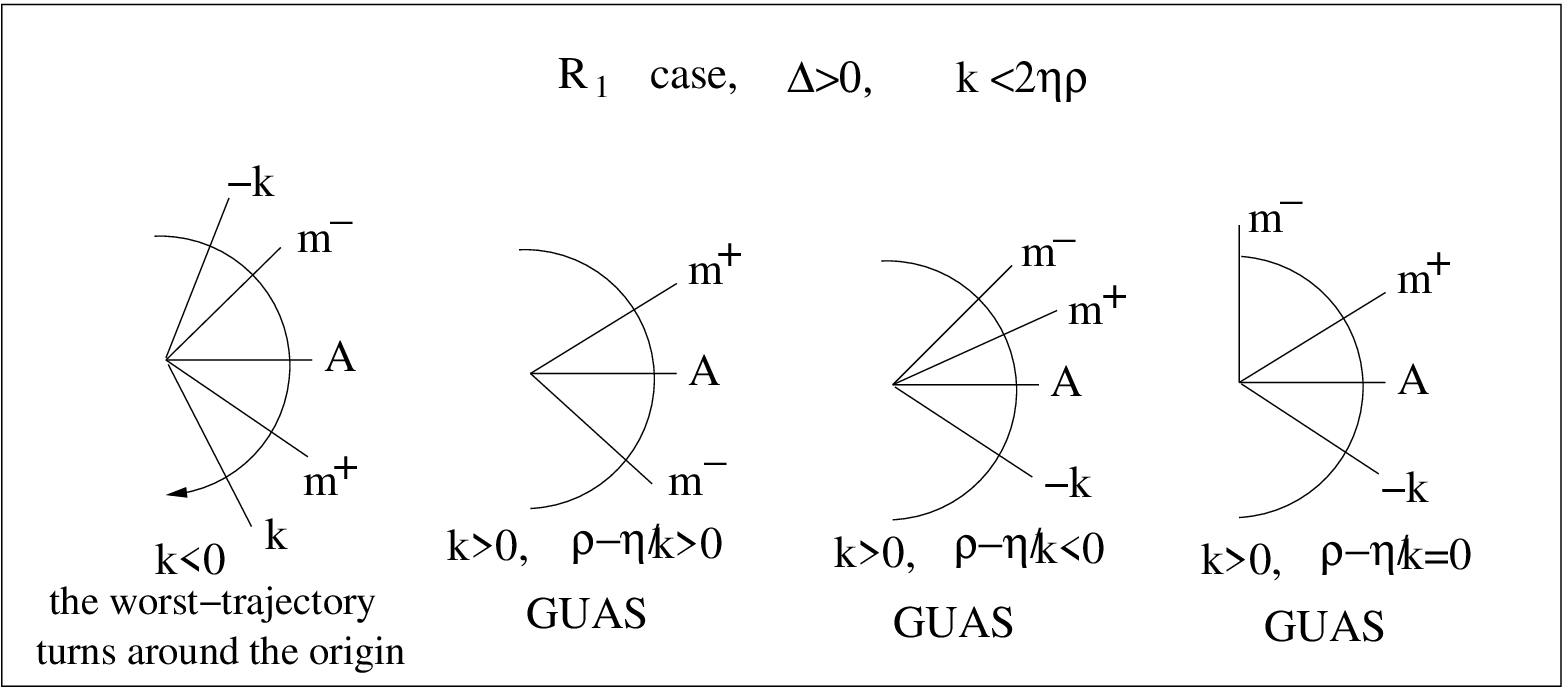}{}{15}
\bd
\item[i)] $k<0$. In this case we have $\chi<0$ and
$\sqrt\Delta>|k|$. It follows that $m^->0>m^+$ and we have:\\\\
%%%%%
{\bf Claim 1.} $k<m^+<m^-<-k.$\\\\
{\bf Proof of Claim 1.}
We have
\bqn
m^+>k \Longleftrightarrow \frac{-k+\sqrt\Delta}{2\chi}>k
\Longleftrightarrow -k+\sqrt\Delta < 2\chi k \Longleftrightarrow
\sqrt\Delta < k(1+2\chi).
\eqn
Since $1+\rho<0$ (cf. Lemma \ref{l-autovalori}), we
have $k(1+2\chi)=k(1+2\rho-2\eta k)=k((1+\rho)+\rho-2\eta k)>0$. Hence
\bqn
m^+>k \Longleftrightarrow
\Delta <
k^2(1+2\chi)^2
\Longleftrightarrow
-4\eta k \chi< 4 \chi k(k+\chi k)
\Longleftrightarrow
-\eta<k+\chi k
\Longleftrightarrow 0<k(\rho+1).
\eqn
The last inequality holds since $\rho+1<0$.
Similarly we prove that
$m^-<-k.$ \quadp
%%%%%%%%%%%%%%%%%%%%%%%%%%%%
\item[ii)] $k>0$ and $\chi>0$. In this case $\sqrt\Delta>|k|.$ It
follows that $m^-<0<m^+$ and we have:\\\\
%%%%%
{\bf Claim 2.} $m^-<-k<m^+<k$.\\\\
{\bf Proof of Claim 2.} We have
$$ m^+>k \Longleftrightarrow \frac{-k+\sqrt\Delta}{2\chi}>k
\Longleftrightarrow -k+\sqrt\Delta > 2\chi k \Longleftrightarrow
\sqrt\Delta> k+ 2 k \chi \Longleftrightarrow \Delta >
k^2(1+2\chi)^2$$
$$\Longleftrightarrow -4\eta k \chi> 4 \chi k(k+
\chi k) \Longleftrightarrow -\eta> k+ k\rho -\eta
\Longleftrightarrow 0> k(\rho+1).$$ Since $(\rho+1)<0,$ the last
inequality is always true. In a similar way we can show that
$m^-<-k.$\quadp
%%%%%%%%%%%%%%%%%%%%%%%%%%%%%%%%%
\item[iii)] $k>0$ and $\chi<0$. In this case $\sqrt\Delta<|k|.$ It
follows that $m^->m^+>0.$\\
\item[iv)] $k>0$ and $\rho-\eta/k=0.$
In this case $m^-=\infty$, $m^+=-\eta>0$.
\ed

In the case {\bf i)}, the worst-trajectory  $\gamma_{x_0}$ rotates
clockwise around the origin (cf. Remark \ref{r-rotation}) and it is a
concatenation of integral curves of $Ax$ (from the line $D^+$ of
equation $x_2= m^+ x_1$ to line $D^-$ of equation $x_2 = m^-
x_1$) and integral curves of $Bx$ otherwise. See again Figure
\ref{f-ruota-R-1}.
Let $t_1$ be the time needed to  an integral curve of $Ax$ to steer a
point from $D^+$ to $D^-$. Recall formula \r{At}. We have
$$
t_1=
\frac{m^+-m^-}{m^+m^-}=\frac{\sqrt{\Delta}}{\eta
k}.
$$
Let $t_2$ be the time needed to the an integral curve of
$Bx$ to steer a point from $D^-$ to $D^+$. Since
$$
\exp(Bt) = \exp(\rho t) \left (
\begin{array}{cc}
\cosh t & \frac{\sinh t}{k}\\
k \sinh t & \cosh t
\end{array}
\right ),$$ we have
$$ t_2 = \mbox{arctanh} \frac{\sqrt{\Delta}}{k(\rho - \frac{\eta}{k})-
\eta}.$$
Defining again ${\cal R}$ as in formula \r{R}, we have similarly to the
previous case
$${\cal R}= \exp(\frac{\sqrt{\Delta}}{k} )\exp(\rho
t_2 )|(1+\frac{m^+\sqrt{\Delta}}{k \eta})(\cosh t_2 +
\frac{m^-}{k}\sinh t_2)|.$$
After simplification we get
$$
{\cal
R}=|(\frac{-k+\sqrt{\Delta}}{-k-\sqrt{\Delta}})\frac{2 k \rho^2 -
Tr(AB)- \sqrt{\Delta}}{2 k \sqrt {\det B}
 (\rho-\frac{\eta}{k})}|\exp(\frac{\sqrt{\Delta}}{k} )\exp(\rho \theta_1),
$$
 i.e.,
$${\cal R}= |(\frac{-k+\sqrt{\Delta}}{-k-\sqrt{\Delta}})\frac{2 k \rho^2 -
    sign(\delta) (Tr(AB)+ \sqrt{\Delta})}{2 k \sqrt {\det B}
 (\rho-sign(\delta)\frac{\eta}{k})}|\exp(\frac{\sqrt{\Delta}}{k} + \rho \theta_1).$$
As in the previous case, the system is GUAS, uniformly stable (but not
GUAS) or unbounded respectively if ${\cal R}<1$, ${\cal R}=1$, ${\cal R}>1$.\\\\
%%%%%%%%%%%%%%%%%%%
Cases {\bf ii)},
{\bf iii)},
{\bf iv)},
are easily studied projecting the system on $\rp$. This is the purpose of
the next section.

%%%%%%%%%%%%%%%%%%%%%%%%%%%%%%%%%%%%%%%%%%%%%%%%%
\subsubsection{The system on the projective space}
In the {\bf R}$_1$-case, it may happen that worst trajectory goes to
the origin without rotating. In this case the system is GUAS. The
description of this fact is contained in the following Lemma, where
we project the system on $\rp$ using its linearity. Recall Lemma \ref{l-w}.
%%%%%%%%%%%%%%%%%%%%%%%%%%%%%%%%%%%%%%%%%%%%%%%%%%%%%%%%%%%
\bl
\label{l-projective}
Assume {\bf (H0)}, $\Delta>0$, $\delta>0$, $k<2\eta\rho$.
Consider the projective space
$\rp$ represented by the semicircle
$\{(x_1,x_2)\in\R^2:~x_1^2+x_2^2=1\mbox{ and }x_1\geq0\}$ with the points
$(0,1)$ and $(0,-1)$ identified.
On $\rp$ represent the eigenvector of $A$ as the point $p_A=(1,0)$, the
eigenvectors of $B$ as the points
\bqn
p_B^+=(\frac{1}{\sqrt{1+k^2}},\frac{k}{\sqrt{1+k^2}} ),~~~
p_B^-=(\frac{1}{\sqrt{1+k^2}},-\frac{k}{\sqrt{1+k^2}} ),
\eqn
and the set $\qq\cap\rp$ as the couple of points
\bqn
p_m^+=(\frac{1}{\sqrt{1+(m^+)^2}},\frac{m^+}{\sqrt{1+(m^+)^2}} ),~~~
p_m^-=(\frac{1}{\sqrt{1+(m^-)^2}},\frac{m^-}{\sqrt{1+(m^-)^2}} ).
\eqn
See Figure \ref{f-proj}.
Let $C_A$ be the connected component to $p_A$ of the set
$\rp\setminus\{p_m^+,p_m^-\}$.
If $p_B^+$ or $p_B^-$ belongs to $C_A$, then the system is GUAS.
\el
\ppotR{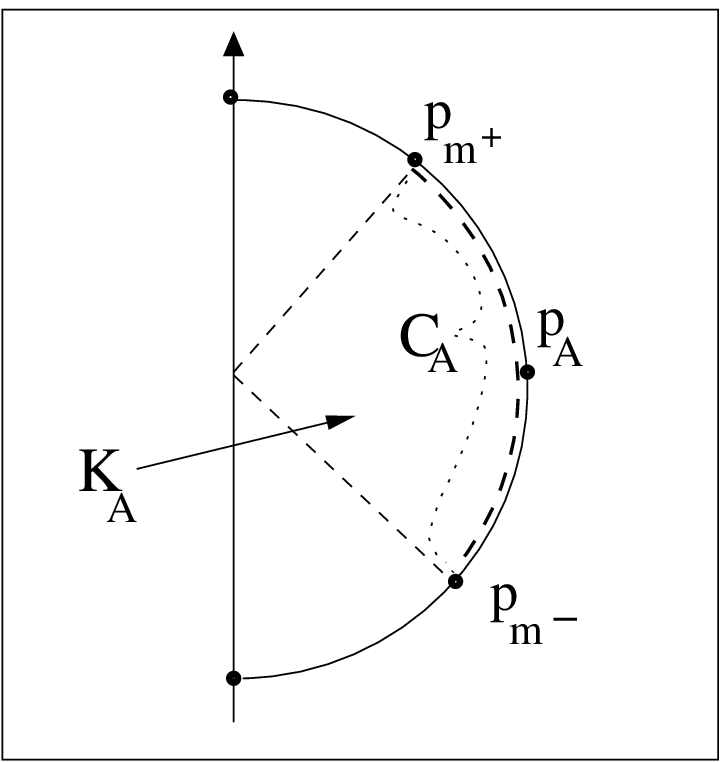}{}{5}
{\bf Sketch of the proof of Lemma \ref{l-projective}.}
By contradiction, if the system is not GUAS, then the worst trajectory  rotates around the
origin and it switches an infinity number of times on the
set $\qq$. Under the hypotheses of the Lemma, once a trajectory enters
the cone
\bqn
K_A=\{x\in\R^2\setminus\{0\}:~~\frac{x}{|x|}\in C_A\}
\eqn
it cannot leave it. Indeed, if in $K_A$ the worst trajectory corresponds
to the vector field $B$ (resp. $A$), it cannot cross the eigenvector of
$B$ (resp $A$) lying in $C_A$, since $p_A$, $p_B^+$ and $p_B^-$ are
stable points for the
system projected on $\rp$). Hence this trajectory cannot rotate around the origin.
Contradiction.
\quadp\\\\
In the cases
{\bf ii)},
{\bf iii)},
{\bf iv)},
the system is GUAS thanks to Lemma \ref{l-projective}.

%%%%%%%%%%%%%%%%%%%%%%%%%%%%%%%%%%%%%%%%%%%%%%%%%%%%%%%%%%%%%%%%%%%%%%%%%%%%%%%
%%%%%%%%%%%%%%%%%%%%%%%%%%%%%%%%%%%%%%%%%%%%%%%%%%%%%%%%%%%%%%%%%%%%%%%%%%%%%%%
%%%%%%%%%%%%%%%%%%%%%%%%%%%%%%%%%%%%%%%%%%%%%%%%%%%%%%%%%%%%%%%%%%%%%%%%%%%%%%%
\subsection{The R$_{0}$ Case}
In this section we study the R$_{0}$-case, i.e., the case in which
the matrix $B$ is also nondiagonalizable. In this case
$\Delta=k^2-4\eta k \rho $. Again we assume
$\Delta>0$, and $k<2\eta\rho$. The set $\qq$ is the union of a
pair of straight lines with angular coefficients given by formula
\r{m00}. Notice that  $\Delta>0$, and $k<2\eta\rho$ implies that
$k<0$ and $m^+<0<m^-$. Indeed if $k>0,$ then
$\Delta= k^2- 4 k \eta \rho = k(k- 2\eta \rho) -2\eta \rho k<0$.
From Lemma \ref{l-autovalori} we know that $B$ has the  unique
eigenvector
$(0,1)$. The relative position of $m^+$, $m^-$ and of the eigenvectors of
$A$ and $B$ is shown in Figure \ref{f-r-0}.
\ppotR{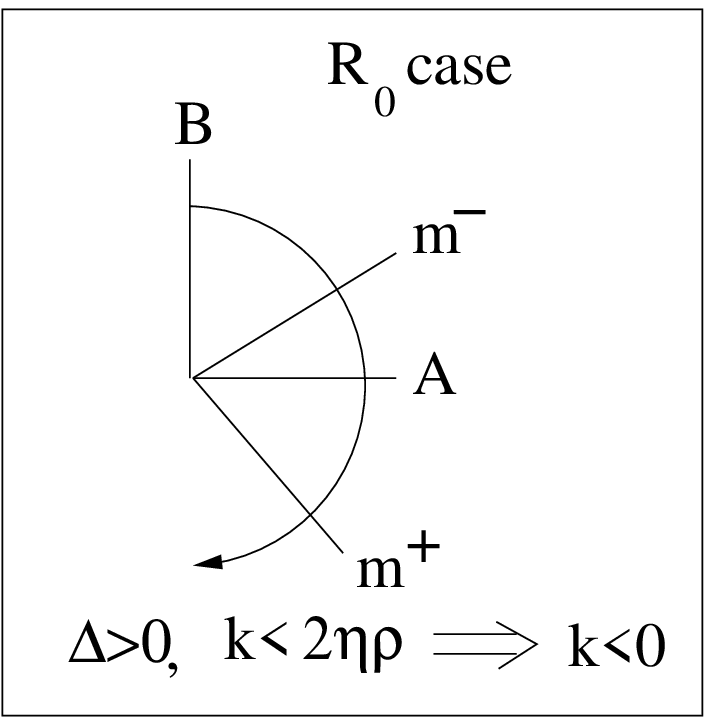}{}{6}
In this case  the worst-trajectory rotates clockwise around
the origin and it is a concatenation of integral curves of $Ax$ from the
line $D^+ =\{(x_1,x_2)\in\R^2\mid x_2 = m^+ x_1\}$ to the line $D^- =\{(x_1,x_2)\in\R^2\mid x_2 = m^-
x_1$\}, and
integral curves of $Bx$ otherwise.
In this case we have
$$
\exp(Bt) = \exp(\rho t) \left (
\begin{array}{cc}
1 &  0\\
k t & 1
\end{array}
\right ),
$$
and the switching times are respectively:
$$t_1= \frac{m^+-m^-}{m^+m^-}=\frac{\sqrt{\Delta}}{\eta k}  \ \ \text{and} \ \
 t_2 = \frac{\sqrt{\Delta}}{k\rho}.$$
Defining again ${\cal R}$ as in formula \r{R}, we have similarly to the
previous cases
\bqn
{\cal R}&=& |\exp(\frac{2\sqrt{\Delta}}{k} )(1+\frac{m^+\sqrt{\Delta}}{k
\eta})|=
|(\frac{-k+\sqrt{\Delta}}{-k-\sqrt{\Delta}})|\exp(\frac{2\sqrt{\Delta}}{k})\nn\\
 &=& |(\frac{-k+\sqrt{\Delta}}{-k-\sqrt{\Delta}})\frac{2 k
\rho^2 - sign(\delta)(Tr(AB)+ \sqrt{\Delta})}{2 k \sqrt {\det B}
 (\rho-sign(\delta)\frac{\eta}{k})}|\exp(\frac{\sqrt{\Delta}}{k}+
\rho \theta_0 ).\nn
\eqn
Again the system is GUAS, uniformly stable (but not GUAS)
or unbounded respectively if
${\cal R}<1$, ${\cal R}=1$, ${\cal R}>1$.
%%%%%%%%%%%%%%%%%%%%%%%%%%%%%%%%%%%%%%%%%%%%%%%%%%%%%%%%%%%%%%%%%%%%%%%%%%%
%%%%%%%%%%%%%%%%%%%%%%%%%%%%%%%%%%%%%%%%%%%%%%%%%%%%%%%%%%%%%%%%%%%%%%%%%%%
%%%%%%%%%%%%%%%%%%%%%%%%%%%%%%%%%%%%%%%%%%%%%%%%%%%%%%%%%%%%%%%%%%%%%%%%%%%
\subsection{The S-case}
In the singular case, the worst-trajectory is always tending to
the origin (see Figure \ref{f-singular}). This is  due to the
fact that the straight line $D^+=\{(h,0)\in\mathbb{R}^2, \
h\in \mathbb{R}\}$ (belonging to $\qq$) coincide with the
eigenvector of $A$ and with
one eigenvector of $B$. Using arguments similar to those of Lemma \ref{l-projective}, one sees that the
worst-trajectory can never
cross $D^+$ and it cannot rotate around the origin.
Hence the system is GUAS. One can also check that in this case
$V(x)=x_1^2+x_2^2$ is a common Lyapunov function.
\ppotR{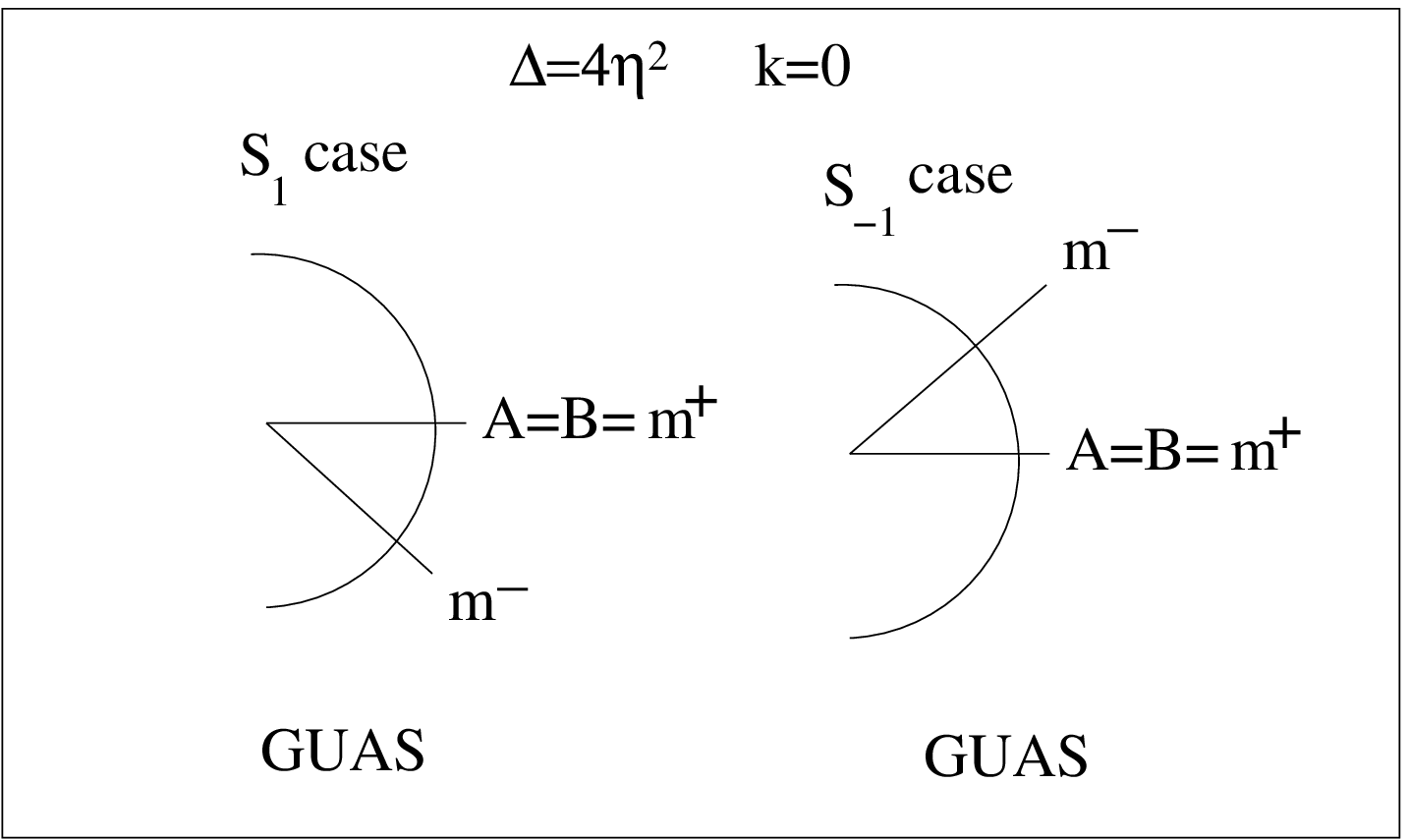}{}{8}

\section*{Acknowledgments}
The authors are grateful to Prof. Pierre Molino, for very helpful
discussions.

\end{document}